\documentclass[12pt]{article}
\title{Henkin constructions of  models with size continuum}
\author{John T. \ Baldwin\thanks{Partially supported
by Simons grant MPS-SCG 418609}\\Department of Mathematics\\University of Illinois, Chicago\and Michael C.\
Laskowski\thanks{Partially supported
by NSF grant DMS-1308546. Both authors acknowledge the support of  DMS
1362974 for visits to Rutgers.}
\\
Department of Mathematics\\University of Maryland}

\usepackage{amssymb}
\usepackage{amsthm}
\usepackage{amsmath}
\usepackage{url}
\usepackage{enumerate}

\def\mcolon{\!:\!}
\def\conc{{\char'136}}
\def\abar{\bar{a}}
\def\bbar{\bar{b}}
\def\cbar{\bar{c}}
\def\dbar{\bar{d}}

\def\ubar{\bar{u}}
\def\vbar{\bar{v}}
\def\wbar{\overline{w}}
\def\xbar{\overline{x}}
\def\ybar{\bar{y}}
\def\zbar{\bar{z}}

\def\phi{\varphi}

\def\D{{\cal D}}

\def\FF{{\bf F}}
\def\H{{\cal H}}

\def\P{{\cal P}}
\def\PP{{\mathbb P}}

\def\T{{\cal T}}

\def\tp{{\rm tp}}

\def\cl{{\rm cl}}
\def\acl{{\rm acl}}
\def\dcl{{\rm dcl}}
\def\pcl{{\rm pcl}}
\def\qcl{{\rm qcl}}

\def\Fa0{{\FF^a_{\aleph_0}}}

\def\bp{{\bf Proof.}\quad}
\def\endproof{\medskip}
\def\<{\langle}
\def\>{\rangle}
\def\o2{{^{\omega} 2}}

\def\n2{{^{n} 2}}

\def\etabar{\overline{\eta}}
\def\sprk{{\rm sprk}}
\newtheorem{Theorem}{Theorem}[subsection]
\newtheorem{Proposition}[Theorem]{Proposition}
\newtheorem{Definition}[Theorem]{Definition}
\newtheorem{Notation}[Theorem]{Notation}

\newtheorem{Remark}[Theorem]{Remark}

\newtheorem{Lemma}[Theorem]{Lemma}
\newtheorem{Corollary}[Theorem]{Corollary}

\newtheorem{Fact}[Theorem]{Fact}

\def\th{Th}
\def\V{\mathbb V}

\newcommand\myrestriction{\mathord\restriction}
\def\mr#1{\myrestriction_{#1}}

\begin{document}

\date{\today}

\maketitle

\section{Introduction}\label{intro}

In the novel {\em White Light} \cite{Rucker}, Rudy Rucker proposes a metaphor for
the continuum hypothesis.  One can reach $\aleph_1$ by a laborious climb up the side of Mt.\ ON, pausing at $\epsilon_0$.
Or one can take Cantor's instantaneous elevator through the center of the mountain.
In this paper, working in ZFC, we take Shelah's elevator, which is a bit slower.
After countably many floors, each with finitely many rooms,
we reach
an object of cardinality $2^{\aleph_0}$.  The underlying construction applies
for finding atomic models, two-cardinal theorems, a collection of continuum many
points that are {\em asymptotically similar} (a weak form of indiscernibility), and a coloring with a Borel square of size continuum.

In his seminal {\em Denumerable models of complete theories},
\cite{Vaught}, Vaught introduced the notion of an atomic
model\footnote{Recall that a formula $\phi(\wbar)$, where $\lg(\wbar)
=n$, is {\em complete} for $T$ if for every formula $\psi(\wbar)$,
$\phi(\wbar)$ {\em decides} $\psi(\wbar)$ in $T$. I.e. $T \vdash
\forall \wbar [\phi(\wbar) \rightarrow \psi(\wbar)]$ or $T \vdash
\forall \wbar [\phi(\wbar) \rightarrow \neg \psi(\wbar)]$. A model $M$
is {\em atomic} if every finite tuple from $A$ satisfies a complete
formula. Here, atomic means $\phi$ is an atom in the Boolean algebra
$F_n(T)$ and has nothing to do with the quantifier rank of the
formula $\phi$.}.  He showed that if the isolated types were
dense\footnote{For every formula $\phi(\xbar)$ consistent with $T$
there is a complete formula $\psi(\xbar)$ such that $T \vdash \forall
\xbar [\psi(\xbar)\rightarrow \phi(\xbar)]$.} in $S(T)$ then $T$ has an
atomic model. Interestingly, \cite{HirShoreSla} show that this
central model theoretic  theorem is not equivalent to any of the
so-called `big five' standard systems of reverse mathematics. Vaught
further showed that a countable atomic model of a complete theory $T$
could be elementarily embedded in every other model; that is, it is
{\em prime}.

The construction of  uncountable atomic models begins with Vaught's
proof \cite{Vaught} that if a countable atomic model has a proper
 atomic elementary extension then it has an atomic elementary
extension of cardinality $\aleph_1$. He constructs a continuous, increasing
sequence of $\omega_1$ countable atomic models and, using the facts
that unions of atomic models are atomic and elementarily equivalent
countable atomic models are isomorphic, deduces the union of the
chain is atomic.  However the construction of atomic models in
cardinals beyond $\aleph_1$ is a long standing problem. The study of
atomic models of complete first order theories translates  to the study of
complete (decides every
$L_{\omega_1,\omega}$-sentence) sentences of $L_{\omega_1,\omega}$
sentences. (See, e.g., Subsection~\ref{reducing} of this paper or Chapter 6 of \cite{Baldwincatmon}.).

Knight \cite{Knightex} showed that construction could stop at
$\aleph_1$; there is a first order theory with no atomic model of
cardinality greater than $\aleph_1$. A series of works
(\cite{Kuekeratomic, LaskowskiShelahatom}) culminating in Hjorth
\cite{Hjorthchar} show that for each countable ordinal $\alpha$ there is
a complete sentence of $L_{\omega_1,\omega}$ that has a model in
$\aleph_\alpha$ but no larger. Thus, it is consistent
that these sentences have no model in the continuum.


Given an atomic model  $M$  of cardinality $\aleph_1$ in a countable
vocabulary, we describe simple sufficient conditions to construct an
elementarily equivalent model $N$ of  cardinality $2^{\aleph_0}$,
which is atomic and Borel.  We modify Henkin's construction  to build a
complete diagram on a family of $2^{\aleph_0}$ variables. The
traditional two steps in a Henkin construction, {\em completeness}, which
ensures that each sentence is decided and {\em Henkin witnesses}, which
ensures that each existential commitment is met, are supplemented by a
 crucial  {\em splitting} stage
which guarantees the final model has the cardinality of the
continuum.

This method generalizes Shelah's construction of  a kind of `tree
indiscernibility', which we call `asymptotic similarity' to give a unified treatment of results in several
areas of model theory.
While we stressed atomic models in the first two paragraphs, the
method applies as well to transfer cardinals in which a type is
omitted and for two cardinal transfers.


We begin by describing the general method in the first five sections.
Section~2 is an overview of both the classical Henkin construction and hints at the new construction.
Section~3 lists a number of desirable properties we might wish the final model satisfied.  Section~4 gives considerably more detail.
There we define finite maximal antichains (fmacs) $A$ of $2^{<\omega}$, $A$-commitments, and generating sequences.
Theorem~\ref{MAIN} of Section~5 is the main result of the paper.

The second half of the paper discusses applications of this technique.  Most of the results are known, but Theorem~\ref{getcontpcl} is new.
Our first application in Subsection~6.1 constructs highly controlled models of theories with {\em trivial definable closure}, which is a notion studied by Ackerman, Freer, and Patel in
\cite{AFP}.  In Subsection~6.2 we introduce the notion of a {\em sufficient pregeometry} and prove, e.g., if $M$ is uncountable and atomic and $(M,\cl)$ is a sufficient pregeometry, then there is an atomic model $N$ of size continuum elementarily equivalent to $M$.  This result immediately entails the new theorem that a pseudominimal theory has an atomic model of size continuum.  In Subsection~6.4 we show that old results of Hrushovski and Shelah from \cite{HruSh334} fit nicely into our rubric.  In particular, if a superstable theory $T$ has an atomic model of size $\aleph_1$, it has an atomic model of size $\beth_1$ (i.e. the continuum).

Section 7 is devoted to streamlining our method under the additional assumption that the theory $T$ has Skolem functions.  In Subsection~7.1 we show that Shelah's
celebrated two-cardinal transfer theorem $(\aleph_\omega,\aleph_0)\rightarrow (2^{\aleph_0},\aleph_0)$
from \cite{Sh37,Sh49} fits this framework.  In Subsection~7.2, we discuss results of Shelah from \cite{Sh522} that describe a cardinal
$\lambda_{\omega_1}(\aleph_0)$ that is large enough so that
any structure $M$ of at least this size can witness arbitrarily long splittings.  As one application, we expound Shelah's proof of the consistency
with $ZFC+2^{\aleph_0}>\aleph_{\omega_1}$ of the statement:
`A sentence of $L_{\omega_1,\omega}$ that has a
model in $\aleph_{\omega_1}$ has one in the continuum.'

This analysis also connects with the philosophical discussion of the
nature of  mathematical explanation. Hafner and Mancosu
\cite{HafnerMancosu} criticized the Resnik and Kushner
\cite{ResnikKushner} assertion that Henkin's proof \cite{Henkin} of
the completeness theorem for first order logic and type theory  is
explanatory. They asked `what the explanatory features of this proof
are supposed to consist of?'. By its explicit  connections with the
deductive system Henkin's original proof was more explanatory of
first order completeness  than G{\"o}del's reduction to propositional
logic \cite{BaldwinChieti}. This paper broadens that debate by noting
that the Henkin construction extends from a transfer from a
syntactic hypothesis  to a semantic conclusion to a transformation
from one model to another.  That is, Henkin's essential contribution
is to explain the ingredients to construct a model. So the
significance of the method is seen in a larger context than the
original proof.

\section{General strategy}\label{gs}
\numberwithin{Theorem}{section} \setcounter{Theorem}{0}

We suppose throughout that we are working with a countable language
$L$ with equality.
Our objective will be to describe  techniques, which are highly
analogous to a Henkin construction of a countable model, for
constructing a model $M$ of size continuum.


%
%
%
%
%
%
%
%
%
%
%
%
%

Classically, the key notion is that of a {\em Henkin set} of formulas, whose definition is rather tedious, but provides the bridge between proof systems and structures.
In their proofs of the completeness theorem both Henkin and G{\" o}del
worked in a framework in which equality was just another relation symbol. And each added an addendum that the proof transferred to the situation where
equality was required to be interpreted as identity. This weakened
 (e.g. Henkin's) conclusion that the model constructed for a vocabulary of
  size $\kappa$ had cardinal $\kappa$ to `at most $\kappa$' or G{\" o}del version for countable languages allowed finite models. Because in our inductive construction we
  will introduce distinct variables that are later forced to be equal,
   we assume predicate logic includes the equality axioms, so
  {\em all witness sets
   will satisfy the usual equality axioms}.  We are not giving a proof of the completeness theorem but transferring the existence of a model with
    specified properties to a model with the same properties but having
    cardinality $2^{\aleph_0}$.

Henkin's most fundamental innovation (e.g.,
\cite{BaldwinChieti}) was to replace
the Skolem functions in G{\"o}del's proof by carefully described constants.
This allowed the transformation from G{\"o}del's universal vocabulary with relation symbols of all
arities to a vocabulary tailored for the topic at hand.

%
%
%
%
%
%

%

%

\begin{Definition}  {\em  Let $L$ be any countable language. Let $Z$ be a distinguished set of indexed variable symbols.
After Henkin,
 $Z$  was viewed as a countably infinite set of constant
symbols.
Here we treat the witnesses as variables so as to encode restrictions on the relations among variables introduced at different levels as transparent validities.


For any $L$-formula $\phi$  with at most $k$ free variables and for
any set of variables $V$, we introduce the notion of a {\em $V$-instantiated
formula}. For any $(v_1,\dots,v_k)\in V^k$, let $\phi(v_1,\dots,v_k)$
be the result of substituting  the variable symbol $v_j$ for the
$j$th free variable for each $j$.
 We call $\phi(v_1,\dots,v_k)$ a $V$-{\em
instantiated formula};
 $Fm(V)$ denotes the set of all  formulas obtained by this procedure.


{\em A witnessed Henkin set} is a subset $\H\subseteq Fm(Z)$ such that:
\begin{itemize}
\item  {\bf  Satisfiable:}  If $\phi(z_1,\dots,z_k)\in\H$, then there is some $L$-structure $N$ and $(a_1,\dots,a_k)\in N^k$ such that $N\models\phi(a_1,\dots,a_k)$.
\item{\bf Completeness:}  For every $\phi\in Fm(Z)$, exactly one of $\phi,\neg\phi\in \H$; and
\item{\bf Henkin witnesses:}  If $\exists w\phi\in Fm(Z)$, then either $\neg\exists w \phi(w)\in \H$ or $\phi(z^*)\in \H$ for some $z^*\in Z$.
\end{itemize}
}
\end{Definition}

It is routine to see that for any witnessed Henkin set $\H\subseteq
Fm(Z)$, the binary relation $z\sim z'$ iff $(z=z')\in\H$ is an
equivalence relation.  As notation, for each $z\in Z$, let $[z]$
denote the image of $z$ under the canonical projection
$\pi:Z\rightarrow Z/\sim$.  The following proposition is proved by a
routine induction on the complexity of formulas; the `Henkin
witnesses' clause is precisely  what is needed to allow quantifiers
to be interpreted correctly.


%

\begin{Proposition}  \label{temp} If $\H\subseteq Fm(Z)$ is a witnessed Henkin set, then there is a unique $L$-structure $M$ with universe $Z/\sim$ that satisfies
$$M\models \phi([z_1],\dots,[z_k]) \qquad \Longleftrightarrow \qquad \phi(z_1,\dots,z_k)\in \H.$$
 In particular, the relation $\sim$ induced by the equality symbol in  $\H$ is a congruence on $Z$.

 Moreover, if $T$ is any
$L$-theory and every
$\phi(z_1,\dots,z_k)\in\H$ is satisfied by some
model $N$ of $T$ (i.e., $N\models \phi(a_1,\dots,a_k)$ for some
$(a_1,\dots,a_k)\in N^k$), then $M$ is a model of $T$.

\end{Proposition}

Note that the whole of the discussion so far does not depend on the size of $Z$!  In the classical construction of a Henkin set,  $Z$ is countably infinite, and
$\H$ is generated by an $\omega$ sequence of formulas $\<\phi_n(\zbar_n):n\in\omega\>$, where, for each $n$, $\zbar_n$ is a subsequence of $\zbar_{n+1}$
and $\phi_{n+1}(\zbar_{n+1})\vdash\phi_n(\zbar_n)$.
{\bf In particular, at each finite stage and for each finite $\zbar\in Z^k$ only `finitely much information' about $\H$ 
is determined.}


%
%

In analogy with this construction, we want to create a template which
can be customized to create a model of size $2^{\aleph_0}$ with
desirable properties. We begin with an indexed set $Z$ of variable
symbols of cardinality $2^{\aleph_0}$, which are subdivided as
$$Z=\bigcup\{Z_s\colon s\ \hbox{a non-empty finite subset of $2^\omega$}\}$$
where each  $Z_s$ is countably infinite and  $Z_t\subseteq Z_s$ whenever $t\subseteq s$.

%
%
%
%
%
%
%
%
%
%

We will construct a witnessed Henkin set $\H\subseteq Fm(Z)$ in
$\omega$ steps. Our subdivision of $Z$ gives rise to sets $Fm(Z_s)$
of instantiated
formulas, whose  intersection with $\H$ yields a
family $\{\H(Z_s)\colon s$ a non-empty finite subset of $2^\omega\}$ of
{\em countable} witnessed Henkin sets.  The restrictions of the congruence
$\sim$ on $Z$
 naturally induce congruences on each $Z_s$.  Thus, exactly as in the classical case outlined above, each of the Henkin sets
$\H(Z_s)$ gives rise to a canonical countable $L$-structure $M(s)$
with universe $Z_s/\sim$. Our construction will ensure that $M(t)$ is
an elementary submodel  of $M(s)$   whenever $t\subseteq s$.

%
%
%
%
%
%

Additionally, the entire Henkin set $\H(Z)$ determines a canonical $L$-structure $M$ with universe $Z/\sim$.
Since any finite tuple $\zbar$ from $Z$ is contained in some $Z_s$, $M$ can be identified with
$$M=\bigcup \{M(s):s \ \hbox{a non-empty finite subset of $2^\omega$}\}$$
In particular, any `finitary information' about $M$ will be inherited from the directed family $\{M(s)\}$ of countable models.
As examples,
\begin{itemize}
\item  $M(s)\preceq M$ for each finite $s\subseteq 2^\omega$, hence for any $T$,  $M\models T$ if and only if some (equivalently, every) $M(s)\models T$;
\item  For $\Delta$ any partial type, $M$ omits $\Delta$ if and only if every $M(s)$ omits $\Delta$; so
\item  $M$ is atomic (Subsection~\ref{atomic}) if and only if every $M(s)$ is atomic.
\end{itemize}

Obviously, if we want to conclude that $M$ has size $2^{\aleph_0}$,
we need some additional mechanism to ensure the {\em construction is
non-degenerate}. In particular, as each $M(s)$ is countable, it would
be very unfortunate if $M(s)=M(t)$ for all finite subsets $s,t$!



 To ensure this, we now introduce
the actual set of variable symbols used in the construction. We will
write $Z=X\cup Y$, where, $X$ is indexed
 as
$\{x_\eta\colon \eta\in 2^\omega\}$ or sometimes we must doubly index
$X$ as $\{x_{\eta,i}\colon \eta\in 2^\omega, i\in\omega\} $.




The intent is that the elements of $X$ are `independent' in some
sense; but at a minimum, we will require that for distinct
$\eta,\eta'$, $x_\eta\neq x_{\eta'}\in\H$\footnote{Or
$x_{\eta,0}\neq x_{\eta',0}\in\H$ in the doubly indexed case.}.
This
will be enough to guarantee that the model $M$ we produce from $\H$ will
have power continuum. The $Y$-symbols are indexed as $\{y_{s,i}\colon
s$ a  non-empty finite subset of $2^\omega$, $i\in\omega\}$ and should be
interpreted as collectively being `material needed to close $X$ into a
model.'   For each non-empty finite subset $s$ of $2^\omega$,
put $X_s:=\{x_\eta\colon \eta\in s\}$ (or
$\{x_{\eta,i},\eta\in s,i\in\omega\}$ in the doubly-indexed case);
put $Y_s:=\{y_{t,i}\colon t\subseteq s, i\in\omega\}$, and $Z_s:=X_s\cup Y_s$.
Visibly, each $Z_s$ is countable and $Z_t\subseteq Z_s$ whenever $t\subseteq s$.


As examples, consider the models $M_{\{\eta\}}$, $M_{\{\eta'\}}$ and $M_s$, where $s=\{\eta,\eta'\}$.
Each of these is a countable, elementary
substructure of $M$.  Thus, in particular, for every constant symbol $c\in L$, there will be natural numbers $i,j$
such that the $Z$-instantiated formulas $y_{\eta,i}=c$ and $y_{\eta',j}=c$ are both in $\H$.  Consequently, $y_{\eta,i}=y_{\eta',j}$ will also be in $\H$,
so $y_{\eta,i}\sim y_{\eta',j}$.  That is, these two variable symbols are identified in both  $M$ and $M(s)$.

  For $s=\{\eta,\eta'\}$, the variables for $M(s)$ are the union of the variables of $M({\eta})$, $M({\eta'})$ and
 $\{y_{s,i}\}$ for $i< \omega$.  The  additional variables $\{y_{s,i}\}$ will close $M(s)$ to be a model.  For example, if we are constructing a
 group,
 then  for some $i\in\omega$, $\H$ would include the $Z$-instantiated formula $x_{\eta}+x_{\eta'}=y_{\{\eta,\eta'\},i}$.
%

%
%


\section{Desirable properties of models}
\numberwithin{Theorem}{subsection} \setcounter{Theorem}{0}

As we are working in a countable language, the existence of
structures, or even models of a consistent first order theory, of size continuum
is not surprising.  Our aim is to identify other desirable properties
of models that do not so obviously have uncountable models but that
can be dovetailed with our construction of a witnessed Henkin set.
Here, we describe some such properties, and the next section will
outline sufficient conditions for a generating sequence and hence a
witnessed Henkin set to admit these properties.

\subsection{Modeling $T$ and omitting types}\label{modot}


We list here the goals of certain conditions on a construction that
will guarantee it yields a model of a given theory $T$ that has the properties we are after. In
Definition~\ref{dense}, we specify how these goals are met in our
situation.

\medskip

\noindent {\bf Modeling $T$:}
As $L$-sentences are themselves $L$-formulas, if we require every
$\phi(\zbar)\in\H$ to be satisfiable in some model of $T$, then the
{\bf Completeness} condition,  each $L$-formula $\phi$ or its
negation is in  $\H$, on a witnessed Henkin set will ensure that
the canonical model $M$ built from $\H$ is a model of $T$.

\medskip

\noindent {\bf Omitting $\Delta$:}
If we want $M$  to omit  a single partial type $\Delta$  we need to require that for any $\zbar\in Z^k$,
there is some  $\delta\in\Delta$ with $\neg\phi(\zbar)\in\H$.  So, if $\H$ is going to be produced in $\omega$ steps,
we need to ensure that every $\zbar\in Z^k$ is `handled' along the way.  Note that, in general, a condition such as `every $\phi\in\H$ is realized in some model
that omits $\Delta$' might not be sufficient to guarantee that $M$ omits $\Delta$.

\medskip

\noindent {\bf Omitting $\{\Delta_m \colon m\in\omega\}$:}
Similarly, if we are given a countable set $\{\Delta_m\}$ of partial
types, in order to ensure that $M$ omits each $\Delta_m$, we need to
ensure that for each pair $(\zbar,m)$, there is a $\delta\in\Delta_m$
for which we enforce that $\neg\delta(\zbar)\in\H$.

\subsection{Atomic models and complete formulas}\label{atomic}

For a complete theory $T$, an $L$-formula $\phi(\xbar)$ is {\em complete with respect to $T$} if:
\begin{itemize}
\item  $T\models \exists\xbar\phi(\xbar)$ and;
\item for every $L$-formula $\delta(\xbar)$, $\phi$ {\em decides}
    $\delta$,
\begin{itemize}
\item
either $T\models \forall\xbar(\phi(\xbar)\rightarrow\delta(\xbar))$;
\item  or $T\models \forall\xbar(\phi(\xbar)\rightarrow\neg\delta(\xbar))$.
\end{itemize}
\end{itemize}
Equivalently, $\phi(\xbar)$ is complete with respect to $T$ if and only if there is a unique complete type extending $\phi(\xbar)$.

A model $M$ of $T$ is {\em atomic} if, for every $n\ge 1$, every  tuple $\abar\in M^n$ realizes a complete formula with respect to $T$.
Not every countable theory $T$ admits an atomic model, but Vaught proved that any two countable, atomic models are isomorphic.  It is easy to see that any elementary submodel of an atomic model is atomic, but the Upward L\"owenheim-Skolem theorem can fail badly -- Hjorth~\cite{Hjorth}
 proved that for any $\alpha<\omega_1$, there are complete theories $T_\alpha$ that have atomic models of size $\aleph_\alpha$, but no larger.
  As it is consistent with $ZFC$ for the continuum to be arbitrarily large in the $\aleph$-hierarchy, we know that we cannot hope to construct an atomic model of size continuum for any of these theories $T_\alpha$. So we must impose some additional hypotheses on $T$ for it to have an atomic model in the continuum.

%

\subsection{$L_{\omega_1,\omega}$-sentences, omitting
types, atomic models}\label{reducing} 

%


We will see that in many cases, the Henkin method will provide
sufficient conditions for building a model of size continuum that is
 atomic, or, in other cases, omits a given countable family of types.  This dual
consequence stems from a fundamental link, discovered independently
by Chang and Lopez-Escobar,  between sentences\footnote{Recall that the logic  $L_{\kappa,\omega}$ allows conjunctions of length less than $\kappa$ but
only finite quantifications; $L_{\infty,\omega} = \cup_{\kappa} L_{\kappa,\omega}$.} $\Phi$ of
$L_{\omega_1,\omega}$ and the omitting of types, which Shelah
extended to atomic models.


%

\begin{quotation}  \label{easytrans}  \noindent Given any sentence $\Phi'$ of $L_{\omega_1,\omega}$ there is a countable language
$L'\supseteq L$, a first-order $L'$-theory $T$, and a partial
$L'$-type $\Delta(w)$ such that the class  of models of $\Phi'$ is
precisely the class of $L$-reducts of models of $T$ that omit
$\Delta(w)$.
\end{quotation}

To see the idea suppose a subformula  $\Phi(\wbar)$ of the sentence
$\Phi'$
 is a countable conjunction of formulas $\phi_i(\wbar)$.
Add a new predicate symbol $R_\Phi(\wbar) $. Let $T$ assert for each
$i$, $\forall \wbar [R_\Phi(\wbar)\rightarrow \phi_i(\wbar)]$ and let
$\Delta(\wbar)$ be the type $\{\neg R_\Phi(\wbar)\} \cup
\{\phi_i(\wbar): i< \omega\}$. Now a model $M$ satisfies $\Phi(\wbar)
\leftrightarrow R_\Phi(\wbar)$ if and only if $M$ omits $\Delta(w)$.
Now hire a secretary who translates the inductive structure of
arbitrary sentence $\Phi'$ into an iteration of extensions of this
sort.

To make the connection with atomic models, we need some further
terminology.

\begin{Definition}  \label{small}  {\em
 An $L_{\omega_1,\omega}$-sentence $\Phi$ is {\em complete} if it has a model and if it decides every $L_{\omega_1,\omega}$-sentence $\Psi$.
 An $L$-structure $M$ is {\em small} if it realizes only countably many distinct $L_{\infty,\omega}$-types over the empty set.
}
\end{Definition}

Recall that each countable model $M$ (in a countable vocabulary) has
a Scott sentence, an $L_{\omega_1,\omega}$-sentence $\Phi_M$, whose
only model is $M$. By the L{\"o}wenheim Skolem theorem $\Phi_M$ is
complete.  Examining the proof of Scott's theorem
(\cite{Keislerbook}) one sees several equivalent statements (see
e.g., Chapter~6 of \cite{Baldwincatmon}):\@
%
%
 an $L_{\omega_1,\omega}$-sentence $\Phi$ is complete if and only
if $\Phi$ is
 $\aleph_0$-categorical if and only if $\Phi$ is a Scott sentence of a
countable $L$-structure. Similar arguments show that an $L$-structure
$M$ is small if and only if it satisfies a complete sentence $\Phi$
if and only if it has a countable $L_{\infty,\omega}$-elementary
substructure if and only if it has a countable
$L_{\omega_1,\omega}$-elementary substructure.

 Shelah
\cite{Sh48} observed:

\begin{Remark}  \label{smalltrans}
If $\Phi$ is a complete $L_{\omega_1,\omega}$-sentence, then there is
a countable language $L'\supseteq L$ and an $L'$-structure
$M'$ such that the class  of models of $\Phi$ is precisely the class of
$L$-reducts of atomic models of $T=\th(M')$.  Conversely, given any complete theory $T$ in a
countable language, there is a complete sentence $\Phi$ of
$L_{\omega_1,\omega}$ whose models are precisely the atomic models of
$T$.
\end{Remark}

\bp  Let $M$ be any countable model of $\Phi$.  For each $k\ge 1$,
define an equivalence relation $\sim_k$ on $M^k$ by $\abar\sim_k
\bbar$ if and only if they have the same $L_{\infty,\omega}$-type
over the empty set.   For each $k$ and $\sim_k$-class $E$, add a new,
$k$-ary predicate symbol $R^k_E$ to $L'$ and let $M'$ be the natural
expansion of $M$, i.e., $M'\models R^k_E(\abar)$ if and only if
$\abar\in E$. Let $T=\th(M')$.

Conversely, given a complete, first order theory $T$, for every $n$ let $\Delta_n(\xbar)$ be the partial type asserting the negation of every complete formula with respect to $T$.
Let $\Phi$ be the $L_{\omega_1,\omega}$-sentence
$$\bigwedge T\wedge\bigwedge_n\forall\xbar\left(\neg\bigwedge \Delta_n(\xbar)\right)$$
 The models of $\Phi$ are precisely the atomic models of $T$.  The completeness of $\Phi$ follows from the uniqueness of countable, atomic models of $T$.
\endproof

%

Because of these observations, the entire subfield of `atomic model
theory' can be considered to be a study of the classes of models of
complete sentences of $L_{\omega_1,\omega}$.  Shelah exploited this
identification by studying atomic models to generalize Morley's
categoricity theorem to $L_{\omega_1,\omega}$ in \cite{Sh87a,Sh87b}.

\subsection{Borel structures}
Following \cite{MontalbanNies}, we say that a structure $M$ is {\em Borel}
if there is a standard Borel space $Z$, a Borel subset $D\subseteq Z$, and a congruence $E\subseteq Z^2$ such
that
\begin{enumerate}
\item  $E$ is a Borel subset of $Z^2$;
\item  The universe of $M$ is $D/E$; and
\item  The pre-image of every  subset of $M^k$ defined by an atomic formula is a Borel subset of $Z$.
\end{enumerate}
If the congruence is the identity, we say that $M$  has an {\em
injective presentation.}

In all cases we consider, the set $Z$ of variable symbols
can be presented as
 a standard Borel space.
As we construct the witnessed Henkin set $\H$ (which yields the entire elementary diagram of $Z$) in $\omega$ steps,
it will follow automatically that the associated model $M$ is a Borel structure, where, moreover $D=Z$.
Typically, however, our methods do not give an injective presentation of $M$.  The one exception to this is in Section~\ref{tttheory},
where we exploit strong hypotheses (trivial definable closure) about the theory that yield an injective presentation.  In that case, we additionally show
that every definable subset of $M^k$ is a finite Boolean combination of open sets.


\subsection{Asymptotic similarity}

 Throughout his career, Saharon Shelah defined
and reaped the benefits from a weakish notion of indiscernibility,
that he used in many varied contexts, including two cardinal
transfer theorems in \cite{Sh37,Sh49}, obtaining perfect squares of
colorings as in \cite{Sh522}, and constructing many models in small,
superstable, non-$\aleph_0$-stable theories. Until now, this notion
was unnamed; we give it a belated baptism as {\em asymptotic
similarity}.

In order to describe this notion we fix some notation for dealing
with sequences from $2^{\omega}$

\begin{Definition}  \label{split} {\em  Fix an integer $\ell$.

\begin{itemize}
\item  A $k$-tuple $(\eta_0,\dots,\eta_{k-1})$ of distinct
    elements from $2^\omega$ {\em splits by $\ell$} if the restrictions
    $\{\eta_i\mr{\ell}\mcolon i<k\}$ to $2^\ell$ are distinct.
\item
Two $k$-tuples $(\eta_0\dots,\eta_{k-1})$
and $(\tau_0,\dots,\tau_{k-1})$ of distinct elements from $2^\omega$ are {\em similar (mod $\ell$)}
if $(\eta_0,\dots,\eta_{k-1})$ splits by $\ell$ and $\eta_i\mr{\ell}=\tau_i\mr{\ell}$ for each $i<k$.
\end{itemize}
}
\end{Definition}

Clearly, every $k$-tuple of distinct elements from $2^\omega$ splits by some $\ell$, and consequently
splits by every $\ell'\ge \ell$; and
similarity (mod $\ell$) is an equivalence relation on the  set of $k$-tuples from $2^\omega$ that split by $\ell$.

\begin{Definition}  \label{asympsim}  {\em
Fix an $L$-structure  $M$.
A subset of $M$, indexed by $\{a_\eta\mcolon \eta\in 2^\omega\}$,
    is {\em asymptotically similar} if, for every $k$-ary
    $L$-formula $\theta$, there is an integer $N_\theta$ such that
    for every $\ell\ge N_\theta$, $$M\models
    \theta(a_{\eta_0},\dots,a_{\eta_{k-1}})\leftrightarrow
 \theta(a_{\tau_0},\dots,a_{\tau_{k-1}})$$
whenever $(\eta_0,\dots,\eta_{k-1})$ and $(\tau_0,\dots,\tau_{k-1})$ are similar (mod $\ell$).
}
\end{Definition}

\begin{Remark}  {\em  Although asymptotic similarity should be thought of as a type of indiscernibility, the indiscernibility is only formula by formula.
For example, consider the structure $M=(2^\omega,U_a)_{a\in 2^{<\omega}}$, where  each $U_a$ is a unary predicate
interpreted as the cone above $a$, i.e., $U_a(M)=\{\eta\in 2^\omega:a\triangleleft \eta\}$.
Then, in $M$, the entire universe $\{\eta:\eta\in 2^\omega\}$ is asymptotically similar, despite the fact that no
two elements have the same 1-type.

This notion of indiscernibles should not be confused with the `tree-indexed
indiscernibles' (which are indiscernible for all formulas in the vocabulary)  in \cite{KKS} which arise from non-superstable
theories and Theorem 3.6 of \cite{Shelahbook}.}
\end{Remark}

\section{Partitions of $Z$ via finite antichains}\label{fmac}
\numberwithin{Theorem}{section} \setcounter{Theorem}{0}


A cursory inspection shows that the set $2^\omega$ is involved in the indexing of elements from $Z$.  We employ the
 standard topology placed on the space $2^\omega $ to describe families of partitions of $Z$.
As notation, for any $a\in 2^{<\omega}$, let $U_a=\{\eta\in 2^\omega:a\triangleleft\eta\}$ and ${\cal U}=\{U_a:a\in 2^{<\omega}\}$.
The {\em standard  topology} on $2^\omega$  is the topology formed by positing that ${\cal U}$ is a base of open sets.

\begin{quotation}  \noindent{\bf Throughout this paper, we will denote elements of $2^{<\omega}$ by lower case roman letters, $a,b,c,\dots$, and we reserve lower case
Greek letters $\eta,\nu,\dots$ for elements of $2^\omega$.}
\end{quotation}

Note that if two elements $a,b\in 2^{<\omega}$ are {\em incomparable,} i.e.,
$a\not\trianglelefteq b$ and $b\not\trianglelefteq a$, then the sets $U_a$ and $U_b$ are disjoint.  A {\em finite, maximal antichain}, abbreviated {\em fmac}
is a finite set $A\subseteq 2^{<\omega}$ in which any two elements are incomparable, and every $b\in 2^{<\omega}$ is comparable to some $a\in A$.
It is easily seen that if $A$ is an fmac, then the sets $\{U_a:a\in A\}$ form a partition of $2^\omega$.  As notation, let $\pi_A:2^\omega\rightarrow A$
denote the projection map, i.e., $\pi_A(\eta)$ is the unique element of $A$ lying below $\eta$.
Curiously, the restriction that $A$ is finite is crucial to obtain a partition of $2^\omega$.
Indeed, if $A$ is any {\em infinite} antichain, then as $2^\omega$ is compact and each of the sets $U_a$ are clopen,  $\{U_a:a\in A\}$ cannot cover $2^\omega$.
Paradigms of fmacs are the sets $2^n$, consisting of all sequences of length $n$, but many other fmacs exist. Our constructions could be done using only the sets $2^n$ but at the cost of suppressing intermediate steps which are fmacs; it is more convenient to do various inductions in the general setting.



{\em We now introduce a second system of variables.}
Given any fmac $A \subseteq 2^{<\omega}$, let $Z_A$ be the following set of variable symbols that are {\em disjoint from $Z$}.
  The indexing on $Z_A$ will parallel that for $Z$.  In particular, $Z_A$ is partitioned into
$X_A\cup Y_A$, $X_A$ is either indexed as $\{x_a:a\in A\}$ or doubly indexed as $\{x_{a,i}:a\in A,i\in\omega\}$,
 and $Y_A=\{y_{t,i}:t\subseteq A,i\in\omega\}$.  For a subset $s\subseteq A$, the sets $X_s$ and $Y_s$ are defined analogously.
Note that in the definition that follows, we build in both the {\bf Satisfiable} condition, as well as a {\em `non-degeneracy'} condition that will imply that the Henkin model we construct has size continuum.

\begin{Definition}  {\em
Let $A \subseteq 2^{<\omega}$ be any fmac.  Define an {\em $A$-commitment} to be a
$Z_A$-instantiated formula

%

$$\phi(\xbar,\ybar),\ \hbox{where $\xbar=\<x_a:a\in A\>$ and $\ybar\subseteq Y_A$}$$
that is  satisfiable in some $L$-structure and with the additional
property that for each $a,a' \in A$, $\phi\vdash x_a\neq x_{a'}$ (or
$x_{a,0}\neq x_{a',0}$ when $X_A$ is doubly indexed). }
\end{Definition}

%
%

To understand the relevance of an $A$-commitment to a Henkin set  $\H$ we are constructing, we need the
 notion of a {\em lifting} $h^*\colon A\rightarrow 2^\omega$ of the fmac $A$ to $2^{\omega}$, which is any (necessarily injective)
mapping satisfying $a\triangleleft h^*(a)$ for every $a\in A$.  Note that any lifting $h^*$ naturally induces an injection, which we also dub $h^*$,
$$h^*\colon Fm(Z_A)\rightarrow Fm(Z)$$
given by replacing each $x_a$ by $x_{h^*(a)}$ and replacing each $y_{s,i}$ by $y_{h^*(s),i}$, where $h^*(s)=\{h(a):a\in s\}$.

Our intent is that if, at some stage of our construction of $\H$ we
include the $A$-commitment $\phi$, we commit ourselves to eventually
making
 $$\{h^*(\phi)\colon\ \hbox{all liftings $h^*:A\rightarrow 2^\omega$}\}$$
 a subset of $ \H$.
 More precisely, we define:

  \begin{Notation}\label{comdef} {\em A {\em commitment} is a pair $(A,\phi)$, where $A$ is an fmac and $\phi$ is an $A$-commitment.
   Each construction will choose a particular set of $A$-commitments (for enough $A$) to
  determine the diagram of $Z$.
  }
\end{Notation}



Given two fmacs $A$ and $B$, we say that {\em $B$ covers $A$},
written $A\le B$, if, for every $a\in A$ there is at least one $b\in
B$ such that $a\trianglelefteq b$. For example, if $n\le m$, then
$2^m$ is a cover of $2^n$.

  If $A\le B$, then a {\em lifting to $B$}
is a (necessarily injective) map $h\colon A\rightarrow B$ satisfying
$a\trianglelefteq h(a)$ for each $a\in A$. Note that if $A\le B$,
then any lifting $h^*\colon A\rightarrow 2^\omega$ factors through
$B$. That is, given any lifting $h^*\colon A\rightarrow 2^\omega$,
 define $h_B\colon A\rightarrow B$ by $h_B(a)=\pi_B(h^*(a))$
(where $\pi_B$ is the  natural projection  from $2^{\omega}$ onto $B$).
Any such $h_B$ is a lifting to $B$, and there is a natural lifting
$h'\colon B\rightarrow 2^\omega$ satisfying $h^*=h'\circ h_B$.

With this in mind, we partially order the set of  commitments by:

\begin{quotation} \noindent  $(A,\phi)\le (B,\psi)$ if and only if $B$ covers $A$ and\footnote{The $\vdash$ means that  $(\forall \zbar) [\psi\rightarrow h(\phi)]$, where $\zbar$ lists the free variables of the formula, is a theorem of the predicate calculus; it is to state this clearly that we work with variables rather than constants.}
$\psi\vdash h(\phi)$ for every lifting $h:A\rightarrow B$.
\end{quotation}

%
%

We say $(B,\psi)$ {\em extends} $(A,\phi)$ when $(A,\phi)\le
(B,\psi)$. Because of our comments about compositions of liftings, it
is evident that whenever  $(B,\psi)$ extends $(A,\phi)$, what $\psi$
commits us to about the $\H$ we will construct is consistent with,
and typically extends what $\phi$ commits us to about $\H$. Thus, if
we have an $\omega$-sequence $\overline {A} =\<(A_n,\phi_n) \colon n\in\omega\>$ of
commitments such that $(A_n,\phi_n)\le (A_{n+1},\phi_{n+1})$ for each
$n$, then let
\begin{quotation}

\noindent $\D_{\overline{A}}:=\{Z$-instantiated formulas $\theta(\zbar)$: for some $n$ (equivalently, for all sufficiently large $n$)
there is some lifting $h^*:A_n\rightarrow 2^\omega$ such that $h^*(\phi_n)\vdash \theta(\zbar)\}$.
\end{quotation}


Visibly, any such set $D_{\overline{A}}$ is closed under logical consequence.
It is
natural to ask for sufficient conditions for a sequence of commitments to determine
a witnessed Henkin set.  More formally:

%

%
%
%

\begin{Definition}  {\em  A {\em generating sequence} is a $\le$-increasing  $\omega$-sequence $\overline {A} = \<(A_n,\phi_n):n\in\omega\>$ of
commitments such that  $ \D_{\overline {A}}$
is a witnessed Henkin set.
}
\end{Definition}

By coupling the discussion in this section with Proposition~\ref{temp}, we see
that if $\overline{A}=\<(A_n,\phi_n):n\in\omega\>$ is a  generating sequence, then
$\D_{\overline{A}}$
uniquely describes a model $M$ of
size $2^{\aleph_0}$.

\section{Sufficient conditions for producing Henkin models of size continuum}\label{suffcond}
\numberwithin{Theorem}{section} \setcounter{Theorem}{0}

We now describe the machinery for constructing a generating sequence.
Even though our construction is in ZFC, cognoscenti will recognize the affinity of our
nomenclature with that of forcing.  We begin by discussing properties of partially ordered sets $(\PP,\le)$ of
commitments.   Note that the `classical Henkin constraints', laid
down in the definition of a witnessed Henkin set, of {\bf
Completeness} and {\bf Henkin witnesses} can be phrased in terms of
showing that certain subsets of $\PP$ are dense and open\footnote{We use Shelah's convention that `more information' puts you `higher up' in $(\PP,\le)$.
Thus, $X$ is {\em dense}
in $(\PP,\le)$ if for every $q \in \PP$, there is an $x \in X$ with
$p \leq x$.  $X$ is {\em open} if $q\in X$ whenever  $q\geq x$ for some $x\in X$.} in $(\PP,\le)$. Additionally, the {\bf Satisfiable}
condition is built into the definition of an $A$-commitment.  The
additional density condition we need to allow us to simultaneously
construct the family $\{M(s)\colon s$ a non-empty finite subset of $2^\omega\}$
of countable models is {\bf Splitting}.


\begin{Definition} \label{splittingdef} {\em
Given any fmac $A$ and any $a\in A$, the {\em splitting of $A$ at
$a$} is
 the fmac $A^{*a}=A\setminus\{a\}\cup\{a\conc 0,a\conc 1\}$.
Clearly, $A^{*a}$  covers $A$, and there are two liftings $h_0,h_1:A\rightarrow A^{*a}$,
distinguished by $h_i(a)=a\conc i$ for $i=0,1$.
Thus,
by the definition of extension, if an $A^{*a}$-commitment $\phi^*$
{\em extends} an $A$-commitment $\phi$
 then $\phi^*\vdash h_0(\phi)\wedge h_1(\phi)\wedge x_{a\conc 0}\neq x_{a\conc 1}$.
 }
\end{Definition}
It is an easy exercise to verify that whenever an fmac $B$ covers $A$, then $B$ can be
obtained by a sequence of splittings at points.
Indeed, the fmac $2^{n+1}$ can be obtained from $2^n$ by a sequence of $2^n$  splittings, one  at each $a\in 2^n$.
The following notation will be used to ensure that appropriate Henkin witnesses are put into a Henkin set.


\begin{Definition}  \label{t(z)}  {\em
Given any fmac $A$ and any finite tuple $\zbar$ from $Z_A$, let $t(\zbar)$ denote the smallest subset of $A$
for which $\zbar\in Z_{t(\zbar)}$.
}
\end{Definition}

Unpacking the definitions, $t(\zbar)$ is the smallest subset of $A$ that satisfies
(1)  If $x_a\in \zbar$, then $a\in t(\zbar)$; and (2) if $y_{s,i}\in \zbar$, then $s\subseteq t(\zbar)$.


\begin{Definition}\label{dense} {\em  A set
$(\PP,\le)$  of commitments, ordered by extension, is {\em
sufficiently dense} if, for every fmac  $A$ and every $A$-commitment $\phi\in
\PP$
we have:

\begin{itemize}
\item{\bf Completeness:}    For every $Z_A$-formula $\psi$, there
    is an $A$-commitment $\phi^*\in\PP$ extending $\phi$ that
    decides $\psi$. By `decides', we mean {\bf either}
    $\phi^*\vdash\psi$ {\bf or} $\phi^*\vdash\neg\psi$;

%
%
%
%
\item{\bf Henkin Witnesses:}  For every $\theta(u,\wbar)$ and
    every $\zbar\in (Z_A)^{\lg(\wbar)}$, there is an
    $A$-commitment $\phi^*\in\PP$ extending $\phi$ such that {\bf
    either} $\phi^*\vdash \forall u\neg\theta(u,\zbar)$ {\bf or}
    $\phi^*\vdash\theta(z^*,\zbar)$ for some $z^*\in Z_{t(\zbar)}$.


\item{\bf Splitting:}  For every $a\in A$ there is an $A^{*a}$-commitment $\phi^*\in\PP$ extending $\phi$.  [In particular,
$\phi^*\vdash h_0(\phi)\wedge h_1(\phi)\wedge x_{a\conc 0}\neq x_{a\conc 1}$.]
\end{itemize}

}

\end{Definition}

Before stating the main theorem, we specify in our context the
properties ensuring the goals laid out at the beginning of
Section~\ref{modot}. They may or may not hold of a particular
$(\PP,\le)$:

\begin{itemize}
\item  {\bf Modeling $T$:}  Given a theory $T$,
if a condition $(A,\phi)\in\PP$, then $\phi$ is satisfiable in
some model of  $T$.

\item  {\bf Omitting
a type
 $\Delta(\wbar)$:}  For every $A$-commitment $\phi\in\PP$ and
    every $\zbar$ from $Z_A$, there is a some $\delta\in\Delta$
    and an $A$-commitment $\phi^*$ extending $\phi$ with
    $\phi^*\vdash\neg\delta(\zbar)$.
\item {\bf Atomic model:}  Given a complete theory $T$, whenever $(A,\phi)\in\PP$, $\phi$ is a complete formula (in its free variables) with respect to $T$.
\end{itemize}

\begin{Theorem} \label{MAIN} Let $T$ be any theory in a countable language.  If there is a sufficiently dense, partially ordered set $(\PP,\le)$ of commitments
 that are each  satisfied in a model of $T$,
then there is a Borel model $M$ of $T$ of size continuum with an asymptotically similar subset $\{a_\eta:\eta\in 2^\omega\}$.  Moreover:
\begin{enumerate}
\item  If $\{\Delta_m\colon m\in\omega\}$ is a countable set of partial
    types\footnote{So the $\Delta_m$ each exemplify a
    $\Delta(\wbar)$ in Definition~\ref{dense}.} and if
    $(\PP,\le)$ satisfies {\bf Omitting $\Delta_m$} for each $m$,
    then such an $M$ can be chosen to omit each $\Delta_m$; and
\item  If $T$ is complete and if $(\PP,\le)$ satisfies the {\bf Atomic model} condition,
then such an $M$ can be chosen to be an atomic model of $T$.

\end{enumerate}
\end{Theorem}

\bp  Fix a distinguished set $Z=X\cup Y$ of variable symbols, for
definiteness\footnote{The $a_\eta$ will be the interpretations of the $x_{\eta,0}$ for $\eta\in 2^\omega $.}, say $X=\{x_{\eta,i}:\eta\in 2^\omega, i\in\omega\}$ and
$Y=\{y_{t,i}:t$ a finite subset of $2^\omega$ and $i\in\omega\}$.

The
following notation will be helpful.  For a fixed $\ell\in\omega$,
consider
 the `standard fmac' $2^\ell$.
In order to consider only finitely many $Y$-variables at each stage,
we distinguish a sufficiently large, finite subset of symbols in
$Z_{(2^\ell)}$.

%
%

Let $$W_\ell:=\{x_{a,i}:a\in 2^\ell, i<\ell\}\cup\{y_{t,i}:t\subseteq
2^\ell, i<\ell\}.$$
%

Note that $W_\ell$ is a finite subset of
$Z_{(2^\ell)}$ and, whenever $\ell\le m$, $h(W_\ell)\subseteq W_m$ for every lifting $h:2^\ell\rightarrow 2^m$.
 We will construct a generating sequence
$\overline{A}=\<(A_n,\phi_n):n\in\omega\>$ from $\PP$  in $\omega$ steps.  We will  dovetail these extensions to obtain the following goals:


\begin{enumerate}[(i)]
\item  All but finitely many  of the `standard fmacs' $2^\ell$ will appear as $A_n$'s in our generating sequence;
\item  To obtain asymptotic similarity, for every formula $\psi(\wbar)$ there is a number $N_\psi$ such that for all
$\ell\ge N_\psi$ there is an $n$ such that $A_n=2^\ell$ and, for every $\zbar$ from $W_\ell$, $\phi_n$ decides
$\psi(\zbar)$;
\item  To show that each of the countable models $M(s)\preceq M$, we require that for every formula $\theta(u,\wbar)$
there is a number $N_\theta$ such that for all $\ell\ge N_\theta$ there is an $n$ such that $A_n=2^\ell$ and, for every $\zbar$ from $W_\ell$,
either $\phi_n\vdash \neg\exists u \theta(u,\zbar)$ or $\phi_n\vdash \theta(y_{t(\zbar),i^*},\zbar)$ for some $i^*\in\omega$ (recall Definition~\ref{t(z)});

\item  Depending on whether we are verifying 1) or 2) there are
    two further conditions.
\begin{enumerate}
\item

For each partial type $\Delta_m(\wbar)$ we are asked to omit,
there will be some $N(m)$ such that for every $\ell\ge N(m)$,
there is an $n$ such that $A_n=2^\ell$ and, for every $\zbar$
from $W_\ell$ (of length $\lg(\wbar)$) there is
$\delta\in\Delta_m$ such that $\phi_n\vdash\neg\delta(\zbar)$;
\item  Finally, if we are asked to produce an atomic model, we require either that every element of $\PP$ be a complete formula, or that for all but finitely many $\ell$,
there is an $n$ such that $A_n=2^\ell$ and, for every $\zbar$ from $W_\ell$, $\phi_n$ entails some complete formula $\eta(\zbar)$.



\end{enumerate}
\end{enumerate}

How can we construct such a generating sequence?  We systematically extend an arbitrary fmac to an $A$ of the form $2^{\ell}$ that satisfies the appropriate condition.  Satisfying (i) is straightforward.  Indeed, given any $(A,\phi)\in\PP$, choose any $\ell$ such that $2^\ell$ covers $A$.
Then, as noted in the discussion above, $2^\ell$ can be obtained from $A$ by a sequence of splittings at points.  So,
it follows from a finite number of applications of {\bf Splitting}
that there is some sequence $\<(B_0,\phi_0),\dots (B_n,\phi_n)\>$ from $\PP$ with $B_0=A$, $B_n=2^\ell$, and
$(B_{i+1},\phi_{i+1})$ extends $(B_i,\phi_i)$ for each $i<n$.


To handle (ii) and (iii),
 fix an enumeration of $L$-formulas
$\{\psi_i(\wbar):i<\omega\}$ and $\{\theta_i(u,\wbar):i<\omega\}$.
For (ii), observe that
as each $W_\ell$ is finite, there are only finitely many
instantiations $\psi_i(\zbar)$ with both  $i<\ell$ and $\zbar$ from
$W_\ell$.  Thus, using the {\bf Completeness} condition on
$(\PP,\le)$ finitely many times,  given any $(A_n,\phi_n)$ with
$A_n=2^\ell$, there is an extension $(2^\ell,\phi_{n+1})\ge
(2^\ell,\phi_{n})$ in which $\phi_{n+1}$ decides every
$\psi_i(\zbar)$ with $i<\ell$ and $\zbar$ from $W_\ell$.

%

Similar remarks concern clause (iii).  Here,  the formulas  $\{\theta_i(u,\zbar):i<\ell\}$ apply, where we
use the {\bf Henkin witnesses} condition finitely often.  Continuing,
again because $W_\ell$ is finite, we can use {\bf Omitting
$\Delta_m$} or {\bf Atomic} to further extend to some
$(2^\ell,\phi_j)\in\PP$ with $j\ge n$ that satisfy iv(a) or iv(b).


Now, once we have handled all of our requirements for the fmac $2^{\ell}$, note that $2^{\ell+1}$ covers $2^{\ell}$, so by finitely many applications of {\bf Splitting} we get
an extension $(A_{n+1},\phi_{n+1})$ with $A_{n+1}=2^{\ell+1}$, thus completing (i) for the next step.
  We repeat the discussion above, but now with the larger $A_{n+1}= B_{\ell+1}$ and a larger (finite) set of formulas $\psi_i(\wbar)\in W_{\ell +1}$
and $\theta_i(u,\wbar)$, for $i < (\ell+1)$.

Continuing this for $\omega$ steps gives us a generating sequence
$\overline{A}=\<(A_n,\phi_n):n\in\omega\>$ from $\PP$.  As cofinally many of the $A_n$'s are $2^\ell$ for increasing $\ell$'s, it follows that $\D_{\overline{A}}$  describes a complete type in the variables $Z$.  The non-degeneracy condition in the definition of a commitment will imply that $\{x_{\eta,0}:\eta\in 2^\omega\}$ are pairwise distinct.  Also, by (ii), this set is easily seen to be asymptotically similar.

%
%


In the construction above, for any witnessed  existential
formula, for all but finitely many $\ell$,  a witness was placed in  $Z_{(2^\ell)}$. Thus, one can check
that if $s$ is a finite subset of $2^\omega$, then $M(s):=\{[z]:z\in
Z_s\}$ is a countable model and $M(s)\preceq M$.
 As well, Clause iv(a) will
imply that $M(s)$ omits each $\Delta_m$, and, in the atomic case,  iv(b) ensures that $M(s)$ is atomic.  As noted in Section~\ref{gs},
knowing that each $M(s)$ omits each $\Delta_m$ or is atomic is enough to conclude that
 $M$ omits each $\Delta_m$
or is atomic.

\section{Applications I - When does an atomic model of size $\aleph_1$ imply one of $\beth_1$?}\label{atapp}
\numberwithin{Theorem}{subsection} \setcounter{Theorem}{0}


In this section, we use the generalized Henkin method to find a number of sufficient conditions on $T$ for which the existence of an atomic model of size $\aleph_1$ implies the existence of an atomic model of size $\beth_1$.
In the first
subsection, we show that if every set is definably closed, a very
 straightforward argument leads from a {\em countable\footnote{Using Theorem~\ref{Vaught},  it is easy to
 see any structure with trivial definable closure is $L_{\omega_1,\omega}$-equivalent to an uncountable structure.}}
model to one in the continuum. In particular, there is no need for the $Y$-variables from our general formulation.
In the second and third subsections we formalize the conditions used in the first in terms of combinatorial geometry and get a general
result which specializes to the goal which motivated this project:
In pseudo-minimal theories \cite{BLSmanymod}, the existence
of an uncountable, atomic model implies one of size continuum.
Then, in the fourth subsection, we move to material that requires much more background and show how the arguments of Hrushovski and Shelah in \cite{HruSh334} can be put into our framework.  There, they prove that if a countable, superstable theory $T$
has an atomic model of size $\aleph_1$, then it has an atomic model of size $\beth_1$.




\subsection{Theories with trivial dcl}   \label{tttheory}

%
%
%
%
%
%
%
%
%

In a series of papers, e.g., \cite{AFP}, Ackerman, Freer, and Patel found that classes of models of
theories with trivial definable closure have some very desirable properties.
Here we note that such theories behave exceptionally well with respect to the Henkin constructions
described in this paper.  In particular, we will see that the {\bf Henkin} and {\bf Splitting} conditions will be
easily satisfied in any model of such a theory.

We begin with a pair of classical definitions.

\begin{Definition}  \label{dclacl}  {\em  Given an $L$-structure $M$ and subset $A\subseteq M$,
an element $b\in M$ is {\em $A$-definable} if there is a formula $\phi(x,\abar)$ with $\abar$ from $A$ for which $b$ is the only solution
in $M$.  The {\em definable closure of $A$}, $\dcl(A)$ is the set of $A$-definable elements of $M$.

Similarly, $b\in M$ is {\em $A$-algebraic} if there is an integer $k$ and a formula $\phi(x,\abar)$ such that $M\models\phi(b,\abar)$
and $M\models\exists^{=k}x\phi(x,\abar)$.  The {\em algebraic closure of $A$}, $\acl(A)$, is the set of $A$-algebraic elements of $M$.
}
\end{Definition}

Clearly, $A\subseteq \dcl(A)\subseteq\acl(A)$ for any subset $A\subseteq M$.
We distinguish structures for which both of these closures are trivial.

\begin{Definition}\label{tottriv}  {\em Fix a countable language $L$.  An $L$-structure $M$
has {\em trivial definable closure} (is $\dcl$-trivial) if
$\dcl(A)= A$ for every subset $A\subseteq M$. }
\end{Definition}

Note that this is very different notion from the usual usage of a trivial closure relation in combinatorial geometry.
Note also  that $\dcl$-triviality is distinct from atomicity. In particular, the theory of countably many independent unary relations is $\dcl$-trivial but has no atomic models.

It is clear that any $\dcl$-trivial structure is infinite, and that
$\dcl$-triviality is a property of the theory of $M$, i.e., if $N$ is
elementarily equivalent to $M$, then $N$ is $\dcl$-trivial if and
only if $M$ is.

The {\em key property} of a $\dcl$-trivial structure $M$ is easy to see:  if $M \models \exists u \phi(u,\cbar) \wedge u \not \in \cbar$, then $\phi(u,\cbar)$ has infinitely many solutions in $M$.  From the key property it is easily seen that $\dcl$-triviality of $M$ is equivalent to $\acl(A)=A$ for every $A\subseteq M$.  In what follows, we will see that $\dcl$-triviality has many equivalent formulations.   A roster of equivalents is given in Fact~\ref{manydcl}.

Constructing models of theories with trivial $\dcl$ is  by far the most straightforward example of our technique, which justifies our considering it first.
The simplicity comes from the fact that we do not require any  $Y$-variables!  But, we must doubly index the x's as $x_{\eta,i}$.

\begin{Definition} \label{nondeg} {\em  Let $N$ be any $L$-structure.  Suppose $\psi(\xbar,\ybar)$ is an $L$-formula with $\lg(\xbar)=k$.  For any $\bbar$ from $N$,
call the definable subset $\psi(N^k,\bbar)$  of $N^k$ {\em non-degenerate} if there exists some $\abar\in \psi(N^k,\bbar)$ with $\{a_1,\dots,a_k\}$ pairwise distinct and disjoint from $\bbar$.
}
\end{Definition}

%
%

\begin{Theorem} \label{tt}  Suppose $M$ is a $\dcl$-trivial structure in a countable language $L$.
There is a model $N$ elementarily equivalent to $M$ of size continuum that satisfies:
\begin{enumerate}
\item  The universe of $N$ is indexed as $2^\omega\times\omega$;
\item  The universe of $N$ can be partitioned as $N=\bigcup_{i\in\omega} A_i$, where, for each $i$,
$A_i=\{a_{\eta,i}:\eta\in 2^\omega\}$ is an asymptotically similar subset;\footnote{In fact, for every finite, strictly increasing sequence $t=(i_1,i_2,\dots,i_k)$ from $\omega$,
the sequences $\{\abar_{\eta,t}:\eta\in 2^\omega\}$ (where $\abar_{\eta,t}=(a_{\eta,i_1},\dots,a_{\eta_{i_k}})$) is
 an asymptotically similar set of $k$-tuples.}
\item  With respect to the natural Polish topology\footnote{The basis consists of  sets of the form $U_a \times \{i\}$where $U_a$ are as in Section~\ref{fmac}.}  on $2^\omega\times\omega$, for every $k$, every
definable subset of $N^k$ is a finite boolean combination of  open sets of $(2^\omega\times\omega)^k$, with the product topology.

\item  If we place the usual  measure\footnote{For any basic open
    $U_{a} \subseteq 2^{\omega}$ with $|a|=n$ let $\mu(U_a)
= \frac{1}{2^n}$ and then extend to $2^\omega\times\omega$ by
 letting $\mu(U_{a} \times \{i\}) = \frac{1}{2^{n+i+1}}$.  In
 fact, if we regard the base set as the locally compact group
 given by pointwise addition on $\omega$ copies of $Z_2^\omega$,
 this is a Haar measure.} on
  $2^\omega\times\omega$, then for every $k$, every
    non-degenerate definable subset of $N^k$ has positive measure
    (with respect to the product measure on
    $(2^\omega\times\omega)^k$.
\item  If, in addition, $M$ is atomic, then we can insist that $N$ be atomic as well;
\item  More generally, if $\{\Delta_m:m\in\omega\}$ is a countable set of types omitted by $M$, then we can insist that $N$ omits each $\Delta_m$ as well.
\end{enumerate}
\end{Theorem}

\begin{Remark}
{\em  In fact, in (3) we can say more  -- the bound on the size of
the boolean combination depends only on $k$, and not on either the
language $L$ or the choice of $L$-structure. That is, there is a
function $k\mapsto n(k)$ with the property that for every countable
$L$ and every $\dcl$-trivial $L$-structure $M$, the associated $N$
has the property that every definable subset of $N^k$ is a boolean
combination of at most $n(k)$ open subsets. }
\end{Remark}

%


{\bf Proof of Theorem~\ref{tt}:}  Fix a $\dcl$-trivial $M$.
We take $Z=X$, where $X$ is doubly
indexed  as $\{x_{\eta,i}:\eta\in 2^\omega, i\in\omega\}$.  To define
our set of commitments, first let $\D_0$ consist of all $L$-formulas
$\phi(\wbar)$ that imply $w_j\neq w_{j'}$ for distinct $j\neq j'$
that are consistent with $T=Th(M)$. For each fmac
 $A$ of $2^{<\omega}$, let $Z_A=\{x_{a,i}\mcolon a\in A,
i\in\omega\}$. Then, for each such $A$, let the set of $A$-commitments $\PP_A$ consist of all
$Z_A$-instantiations of  formulas
 $\phi(\wbar)\in\D_0$ by a tuple $\zbar$ of
{\em distinct} elements of $Z_A$.

Let $(\PP,\le)$
be the poset with universe $\PP=\bigcup\{\PP_A\colon \
\hbox{$A$ an fmac of $2^{<\omega}$}\}$ and where
 $\le$ is the  extension
relation from Section~\ref{suffcond}. We show that {\bf Completeness},  {\bf Henkin witnesses},
and {\bf Splitting} conditions follow easily:
  Fix any fmac
$A$ and any $A$-commitment\footnote{We sometimes abuse notation by identifying $\PP_A$ with the formulas that occur as second coordinates of the pairs.}
$\phi(\xbar)\in\PP_A$. As $\phi(\xbar)$
is consistent with $Th(M)$, choose $\cbar$ from $M$ such that
$M\models\phi(\cbar)$.

\medskip
\noindent
{\bf Completeness:}  Given a $\psi(\zbar)$, where $\zbar$ is a subsequence of $\xbar$, we will show it is decided.  Let $\bbar$ be the corresponding subsequence of $\cbar$.
Now, if $M\models \psi(\bbar)$, then put $\phi^*:=\phi(\xbar)\wedge\psi(\zbar)$, and put $\phi^*:=\phi(\xbar)\wedge\neg\psi(\zbar)$ otherwise.
%
%
%
\endproof

\medskip
\noindent
{\bf Henkin witnesses:}  We must satisfy the condition for an arbitrary $\theta(w,\zbar)$ with $\zbar$ a subsequence of $\xbar$.   Let $t:=t(\zbar)$ be the set of $a \in A$ such that for some $i$, a variable $x_{a,i}$ appears in $\zbar$.
As above, let $\bbar$ be the subsequence of $\cbar$ associated to $\zbar$.   There are three cases.  First, if $M\models\neg\exists w \theta(w,\cbar)$,
then, put $\phi^*:=\phi(\xbar\ybar)\wedge\neg\exists w\theta(w,\zbar)$.  Then $\cbar$ witnesses that $\phi^*$ is an $A$-commitment and it is evident that $(A,\phi^*)$ extends $(A,\phi)$.

Second, suppose $M\models \theta(c,\cbar)$ for some $c\in\cbar$.  Let $z^*$ be the (unique) element of $\zbar$ corresponding to $c$.  Then
$\phi^*:=\phi(\xbar)\wedge\theta(z^*\zbar)$ is in $\PP_A$ and extends $\phi(\xbar)$.

Finally, suppose $M\models \exists u\theta(u,\cbar)\wedge\bigwedge u\not\in\cbar$.  Then, by the key property of $\dcl$-triviality, choose $b^*\in M\setminus\cbar$ such that $M\models\theta(b^*,\cbar)$.
Choose any $a\in t(\zbar)$ and $j\in\omega$ such that $x_{a,j}\not\in\xbar$ and put
$$\phi(x_{a,j}\xbar):=\phi(\xbar)\wedge\psi(x_{a,j}\zbar)\wedge\bigwedge x_{a,j}\not\in\xbar$$
Then $b^*\cbar$ witnesses that $\phi^*\in\PP_A$, which visibly extends $\phi$.

\medskip
\noindent
{\bf Splitting:}
Choose any $a\in A$.
To handle this case, we start with a Claim, whose proof is an easy induction on $k$; the key property yields the case  $k=1$:

\medskip\noindent
{\bf Claim.}  For every $k\ge 1$, for every $\phi(\xbar)\in\D_0$, and for every partitioning of $\xbar=\ubar\vbar$ with
$\lg(\ubar)=k$, then for every $\bbar$ from $M$ such that $M\models\exists\ubar \phi(\ubar,\bbar)$, there is an infinite, pairwise disjoint set
$\{\cbar_j:j\in\omega\}\subseteq M^k$ of realizations of $\phi(\ubar,\bbar)$.

Given the Claim, partition the variables of $\phi(\xbar)$ into two disjoint subsequences $\xbar=\xbar_a\xbar^*$,
where $\xbar_a$ consists of all $x_{a,i}\in\xbar$, while $\xbar^*$ consists of all $x_{a',i}\in\xbar$ with $a'\neq a$.
This partition induces a partition of our realizing sequence $\cbar$ into $\cbar_a\bbar$, where $\cbar_a$ corresponds to $\xbar_a$, while
$\bbar$ corresponds to $\xbar^*$.
Put $$\phi^*(\xbar_{a\conc 0},\xbar_{a\conc 1},\xbar^*):=\phi(\xbar_{a\conc 0},\xbar^*)\wedge\phi(\xbar_{a\conc 1},\xbar^*)\wedge\ \hbox{`$\xbar_{a\conc 0},\xbar_{a\conc 1},\xbar^*$ are distinct'}$$
Then the Claim implies that $({A^{*a}},\phi^*) \in\PP_{A^{*a}}$, and is as required.
\endproof


Now, with our density conditions satisfied, the existence of a model $N$ follows from Theorem~\ref{MAIN}.
By our choice of $\D_0$, the congruence $\sim$ on $Z=X$ is trivial, which establishes Clause~1) and the partition of Clause~2).
The remaining Clauses are established by the properties  guaranteed by Theorem~\ref{MAIN} and the footnotes.
\endproof

\subsection{Sufficient pregeometries}

 In this and the following
subsection we study the effect of having an atomic model  that is equipped with a well behaved closure relation.
In this subsection we give a sufficient  set of conditions on a closure relation of an atomic model $(M,\cl)$ to
allow for the construction of an elementarily equivalent atomic model of size continuum.
As an application, in the next subsection we
prove a new result: among pseudo-minimal theories, the existence
of an uncountable, atomic model implies one of size continuum.

Although we have cast our results in terms of the existence of atomic models, they translate to  complete sentence of $L_{\omega_1,\omega}$ as in Section~\ref{reducing}
(equivalently for countable, first order theories that omit a given type).


\begin{Definition}  \label{fla} {\em
 Let $M$ be any $L$-structure.  A {\em formula-based closure relation on $M$} is a function $\cl:\P(M)\rightarrow\P(M)$ satisfying for all $A,B\subseteq M$,
 $A\subseteq \cl(A)$;
$A\subset B$ implies
$\cl(A)\subseteq\cl(B)$; $\cl(\cl(A))=\cl(A)$; and whenever $a\in \cl(B)$, then there is a finite tuple $\bbar$ from $B$ and a formula $\phi(x,\ybar)\in\tp(a\bbar)$
such that $a'\in\cl(\bbar')$ whenever $M\models \phi(a',\bbar')$.
}
\end{Definition}

Formula-based closure relations abound in model theory.  Examples include {\em equality} $(M,=)$, where $\cl(A)=A$ for all $A\subseteq M$, {\em definable closure}
$(M,\dcl)$, and {\em algebraic closure}  $(M,\acl)$.  Additionally, in the next subsection we introduce {\em pseudo-algebraic closure} $(M,\pcl)$, which is well behaved
whenever $M$ is atomic.  In order to apply our methods, we need our formula-based closure relation to satisfy more properties.


\begin{Definition} \label{sufficient}
{\em Consider a formula-based closure relation $(M,\cl)$ on an arbitrary infinite $L$-structure.
We call $(M,\cl)$ {\em sufficient} if the following additional conditions hold:
\begin{enumerate}
\item  `Exchange:' i.e., if $a\in\cl(Bc)\setminus \cl(B)$, then $c\in \cl(Ba)$;
\item `Extendible\footnote{If any of $\dcl$, $\acl$,  or $\pcl$ are not extendible, the Scott sentence of $M$ has exactly one model.}:' 
There is $a\in M\setminus\cl(\emptyset)$; and
\item  `Weak homogeneity:'  For all finite $\bbar$ and $L$-formulas $\phi(w,\bbar)$, if there is  $a\not\in\cl(\bbar)$ with $M\models\phi(a,\bbar)$,
then for every finite
 $E\subseteq M$, there is $a'\not\in\cl(E)$ that also satisfies $M\models\phi(a   ',\bbar)$.
\end{enumerate}
}
\end{Definition}



A closure relation that satisfies Exchange is also known as a {\em pregeometry} or a {\em matroid}.  It is well known that pregeometries give rise to a well behaved notion of dimension.  In particular, for any set $B$, any two maximal independent subsets of $\cl(B)$ have the same cardinality.  One of many introductions to the role of

\begin{Remark}   \label{quote}  {\em We say $\abar$ is {\em independent over $E$} if for every $i< \lg (\abar)$, $a_i \not \in \cl(\abar -\{a_i\} \cup E)$.  A routine induction shows that the `Weak homogeneity' condition implies that for every $n$, every $\psi(\wbar,\bbar)$, if there is
an $n$-tuple $\abar$ independent over $\bbar$ with $M\models\psi(\abar,\bbar)$, then for every finite $E$, there is $\abar'$ independent over $E$
with $M\models\psi(\abar',\bbar)$.  Also, coupled with `Extendible', 
we conclude that $M$ contains an infinite independent subset $I$.
Moreover, for any $L$-formula $\phi(w,\bbar)$, either $\phi(M,\bbar)\subseteq\cl(\bbar)$, or for every finite set $E$, $\phi(M,\bbar)$ contains an infinite, $E$-independent subset.
}
\end{Remark}

%


Examples of sufficient pregeometries are common.  A structure $(M,=)$ has a sufficient pregeometry if and only if $M$ has trivial $\dcl$.
 If $T$ is strongly minimal, weakly minimal, o-minimal, or has SU-rank 1, then $(M,\acl)$ is a pregeometry for any model of $T$.  Moreover, an easy compactness argument shows that any (infinite) model $M$ of such a theory has a proper, elementary extension $N$ for which $(N,\acl)$ is sufficient.
In the next subsection we prove that whenever a pseudo-minimal theory has an uncountable atomic model, then $(M,\pcl)$ is sufficient for every atomic model.
For now, we content ourselves with the following result.

\begin{Theorem}  \label{big}  Suppose $(M,\cl)$ is a sufficient pregeometry.  Then there is a Borel model $N\equiv M$ of size continuum
with   a $\cl$-independent, asymptotically similar subset $\{a_\eta\mcolon \eta\in \o2\}$ from $N$.
Moreover, if $M$ is atomic (with respect to $Th(M)$) then we may additionally choose $N$ to be  atomic.
More generally, if $\{\Delta_m(\wbar_m):m\in\omega\}$ is a countable set of partial types, each of which is omitted in $M$, then we may additionally require that
$N$ omits every $\Delta_m$.
\end{Theorem}

\bp  In this application, it is helpful to doubly index the $X$-variables.  That is, take as variables $X=\{x_{\eta,i}:\eta\in 2^\omega,i\in\omega\}$,  as usual, $Y=\{y_{s,i}:
s\subseteq 2^\omega$ finite, $i\in\omega\}$ and $Z=X\cup Y$.  The double indexing of the $X$-variables is needed since a typical model (e.g., some $M_{\eta}$)
may have an infinite, independent subset.

As {\bf notation}   , for any fmac $A$,  any non-empty subset $t\subseteq A$, and any $\xbar\in X_A$,
 $\xbar_t$ denotes the subsequence of $\xbar$ from $X_t$, i.e., an element $x_{a,i}\in\xbar$ is an element of $\xbar_t$ if and only if $a\in t$.
Similarly,  for any $\ybar\in Y_A$, $\ybar_t$ is the subsequence of $\ybar$  from $Y_t$, i.e.,  for $y_{s,i}\in \ybar$, $y_{s,i}\in\ybar_t$ if and only if $s\subseteq t$.

For any fmac $A$, let $\PP_A$ denote all $Z_A$-instantiated formulas $\phi(\xbar,\ybar)$ where $\xbar\in X_A$, $\ybar\in Y_A$ and there are sequences $\cbar,\bbar$ from $M$ satisfying:
\begin{enumerate}
\item $M\models \phi(\cbar,\bbar)$;
\item  $\cbar$ is $\cl$-independent; and
\item For each $t\subseteq A$,
$M\models \forall\xbar\forall\ybar(\phi(\xbar,\ybar)\rightarrow \ybar_t\subseteq\cl(\xbar_t))$ (cf., `Formula-basedness')
\end{enumerate}

As usual, let $(\PP, \le)$ be the poset with universe
$$\PP=\{(A,\phi):A\ \hbox{is a fmac and $\phi\in\PP_A\}$}$$
and $\le$ is the usual extension relation.
We argue that $(\PP,\le)$ satisfies   {\bf Completeness}, {\bf Henkin witnesses}, and {\bf Splitting.}

\medskip

Fix an fmac $A$ and an $A$-commitment $(A,\phi(\xbar,\ybar))\in \PP_A$.   Choose finite tuples $\cbar,\bbar$ from $M$ witnessing that $\phi\in\PP_A$.

\medskip
\noindent{\bf Completeness:}
Choose any $\psi(\zbar)$ with $\zbar$ from $Z_A$, which we may assume is a subsequence of $\xbar\ybar$.  Let $\dbar$ be the corresponding subsequence of $\cbar\bbar$.
There are now two cases:  If $M\models\psi(\dbar)$, then put $\phi^*(\xbar\ybar):=\phi(\xbar\ybar)\wedge\psi(\zbar)$; and put
$\phi^*(\xbar\ybar):=\phi(\xbar\ybar)\wedge\neg\psi(\zbar)$ otherwise.  In either case, the same pair $\abar\bbar$ demonstrate that $\phi^*\in\PP_A$.
\endproof

\medskip
\noindent
{\bf Henkin witnesses:}  Choose $\theta(w,\zbar)$ with $\zbar$ from $Z_A$, which we may again assume  is a subsequence of $\xbar\ybar$.
As above, let $\dbar$ be the subsequence of $\cbar\bbar$  corresponding to $\zbar$, and
in the notation of Definition~\ref{t(z)} as amplified just above, let $t=t(\zbar)\subseteq A$.
 There are now three cases.  First, if $M\models\neg\exists w \theta(w,\dbar)$,
then put $\phi^*(\xbar\ybar):=\phi(\xbar\ybar)\wedge\neg\exists w\theta(w,\zbar)$.

Second, suppose there is $h\in \cl(\cbar_t)$ such that $M\models\theta(h,\dbar)$.
  By  `formula-basedness' choose
a formula $\delta(w,\xbar_t)\in\tp(h,\cbar_t)$ such that any realization of $\delta(w,\cbar_t)$ in $M$ implies
$w\in\cl(\cbar_t)$.  Choose $i$ such that $y_{t,i}\not\in\ybar$.
Put $$\phi^*(\xbar,\ybar y_{t,i}):=\phi(\xbar,\ybar)\wedge\theta(y_{t,i},\zbar)\wedge\delta(y_{t,i},\xbar_t)$$
That $\phi^*\in\PP_A$ is witnessed by appending $h$ to $\bbar_t$.

Third, suppose there is $h\in M\setminus \cl(\cbar_t)$ such that $M\models\theta(h,\dbar)$.
Then, clearly, $\{h\}\cup\cbar_t$ is independent.  Choose any  $i\in\omega$ such that $x_{t,i}\not\in\xbar_t$.
Put $$\phi^*(x_{t,i}\xbar,\ybar):=\phi(\xbar,\ybar)\wedge\theta(x_{t,i},\zbar)$$
By Weak Homogeneity choose $c^*\not\in \cl(\cbar\bbar)$ with $M\models\theta(c^*,\dbar)$.
As $\phi^*$ is witnessed by $c^*\abar\bbar$, it follows that $\phi^*\in\PP_A$ and extends $\phi$.
\endproof

\medskip
\noindent
{\bf Splitting:}
Choose any $a\in A$ and let $A^-=A\setminus\{a\}$.
 Partition the variables of $\zbar=\xbar\ybar$ into four disjoint subsequences:
\begin{itemize}
\item $\xbar_a$ is the subsequence of $\xbar$ consisting of all $x_{a,i}\in\xbar$;
\item $\xbar_0$ is the subsequence of $\xbar$ consisting of all $x\in X_{A^-}$;
\item  $\ybar_a$ is the subsequence of $\ybar$ consisting of all $y_{s,i}\in\ybar$ for which $a\in s$; and
\item  $\ybar_0$ is the subsequence of $\ybar$ consisting of all $z\in Z_{A^-}$ (i.e., whose coordinates do not mention $a$).
\end{itemize}

As notation, let $\cbar_a,\cbar_0,\bbar_a,\bbar_0$ denote the subsequences of $\cbar\bbar$ corresponding to $\xbar_a,\xbar_0,\ybar_a,\ybar_0$, respectively.
Put $\psi(\xbar_a,\xbar_0,\ybar_0 ):=\exists \ybar_a \phi$.  Then $M\models\psi(\cbar_a,\bbar_0, \cbar_0)$ as witnessed by $\bbar_a$.  Furthermore, $\cbar_a\cbar_0$ form a partition of $\cbar$ and hence are independent.  Thus, by condition 3)  $\bbar_a\bbar_0\subseteq \cl(\cbar_a\cbar_0)$ and $\cbar_a$ is independent over $\cl(\bbar_0\cbar_0)$.

By Remark~\ref{quote} choose $\cbar_a'$ from $M$ realizing $\psi(\xbar_a,\bbar_0\cbar_0)$ and independent from all of $\cbar\bbar$.  In particular, $\cbar_a'$ is disjoint from $\cbar_a$.
By choice of $\psi$, choose $\bbar_a'$ from $M$ such that $M\models\phi(\cbar_a',\cbar_0,\bbar_a',\bbar_0)$.  It follows that $\bbar_a'\subseteq\cl(\cbar_a'\cbar_0\bbar_0)$.  It is easily checked that these tuples witness: $[h_0(\phi)\wedge h_1(\phi)\wedge \xbar_{a\conc 0}\cap \xbar_{a\conc 1}=\emptyset]\in \PP_{A^{*a}}$.
\endproof

\subsection{Pseudominimal Theories}\label{pcl}

In a series of papers, the authors and Shelah have attempted to determine whether every $\aleph_1$-categorical, complete sentence $\Phi$ of $L_{\omega_1,\omega}$
has a model of size continuum.   By the reductions in Subsection~\ref{reducing}, this is equivalent to asking whether a complete first order theory $T$ that has a unique atomic model of size $\aleph_1$ must also have an atomic model of size continuum.

To analyze this problem, in \cite{BLSmanymod}, we
 introduced a new notion of closure,  which we dubbed  {\em pseudo-closure},
 shortening pseudo-algebraic closure,
that is appropriate for the study of atomic models of a first order
theory. We proved that if pseudo-closure  fails exchange in a strong
way on the class of atomic models of a theory $T$ then $T$ has
$2^{\aleph_1}$ atomic models of cardinality $\aleph_1$. We give a
slightly simplified account of pseudo-minimality which is adequate for the
applications. Here we show that if $T$ has an uncountable atomic
model that is pseudo-minimal, then there is an atomic model of $T$ in
the continuum.

%
%
%


%
%
%
%
%


\begin{Definition}\label{psalgdef} {\em Let $M$ be an atomic model and suppose $a,\bbar$ are  from $M$.
 We say $a$ is {\em pseudo-algebraic over $\bbar$ in $M$}, written $a\in\pcl(\bbar)$, if every elementary substructure $N\preceq M$ that contains $\bbar$ also contains $a$.
}
\end{Definition}

  We showed in \cite{BLSmanymod} that pseudo-algebraicity in atomic models is formula-based and a property of the theory as opposed to a particular model.
That is, if $M$ and $M'$ are elementarily equivalent atomic models, $\abar,\bbar$ and $\abar',\bbar'$ are from $M$ and $M'$, respectively, whose pairs realize the same complete formula,
then $\abar\in\pcl(\bbar)$ in $M$ if and only if $\abar'\in\pcl(\bbar)$ in $M'$.
Also, Lemma~2.6 of \cite{BLSmanymod} implies that if $M$ is atomic, then $(M,\pcl)$ satisfies  the `Weak homogeneity' clause from Definition~\ref{sufficient}.

 Using this notion we can immediately add a clause to an old theorem of Vaught.

 \begin{Lemma}  \label{Vaught}  Let $T$ be a complete theory in a countable language that has an atomic model.  The following notions are equivalent:
 \begin{itemize}
 \item $T$ has an uncountable atomic model;
 \item   the countable atomic model has a proper atomic
     extension;
 \item  the countable atomic model is not minimal; and the new
 \item  $\pcl(\emptyset)\neq M$ for some/every atomic model.
 \end{itemize}
 \end{Lemma}

 \begin{Definition}
{\em Let $M$ be an atomic model and suppose $T$ satisfies the
conditions of Lemma~\ref{Vaught}. We say  that $T$ is {\em pseudominimal}
if $(M,\pcl)$ satisfies Exchange for some/every atomic model $M$ of $T$.
That is, for every finite set $C$
from $M$ and elements $a,b\in M$, if $b\in\pcl(Ca)$ but
$b\not\in\pcl(C)$, then $a\in\pcl(Cb)$.
 }
 \end{Definition}

Thus, a complete theory $T$ satisfying the hypotheses of Lemma~\ref{Vaught} is pseudominimal if and only if
$(M,\pcl)$ is a sufficient pregeometry for some/every atomic model $M$ of $T$.
%

The following  new Theorem is a culmination of our previous results.
It follows immediately from Lemma~\ref{Vaught}, the note above, and Theorem~\ref{big}.

 \begin{Theorem}\label{getcontpcl}
If  a countable first order theory $T$ has an atomic pseudominimal
model $M$ of cardinality $\aleph_1$ then there is an atomic pseudominimal model
$N$ of $T$ with cardinality $2^{\aleph_0}$.

Equivalently, if the models of a complete sentence $\Phi$ in
$L_{\omega_1,\omega}$ are pseudominimal and $\Phi$ has an uncountable
model, it has a model in the continuum.
\end{Theorem}

%
%
%
%
%
%


\noindent
Whereas Theorem~\ref{getcontpcl} is of general interest, we note a special case.  It is an easy exercise to
prove that any weakly minimal theory $T$ with an uncountable atomic model is pseudominimal.  Thus, Theorem~\ref{getcontpcl}
gives a  proof that such a theory has an atomic model of size continuum (a second is Theorem~\ref{ssthm}).

As an example of pseudominimality, Zilber \cite{Zilbercatex, Baldwincatmon} introduced the abstract notion of a
quasiminimal (excellent) class and proved such classes are categorical in all
uncountable powers.  In general, these classes are axiomatized in $L_{\omega_1,\omega}(Q)$ (\cite{Kirbyqm}) and the quasiminimal closure is distinct from
our notion of $\pcl$.  However, in some cases, most notably \cite{ZilberBays},
the study of covers of certain algebraic groups  e.g. \cite{ZilberBays, Baysthesis},  the countability of the quasiminal closure is expressible in
$L_{\omega_1,\omega}$ and then $\pcl = \qcl$.

\subsection{Stable and superstable theories} \label{ss}

Stable theories give rise to a well-behaved notion of independence,
namely non-forking. Using this tool in conjunction with the methods
of this paper,  Hrushovski and Shelah \cite{HruSh334} obtain the
following transfer theorem:

\begin{Theorem} \label{ssthm} Suppose $N$ is an uncountable model of a superstable theory $T$ in a countable language.  Then there is an atomic model $M$ of $T$ of size
continuum that has an asymptotically similar subset $\{a_\eta:\eta\in 2^\omega\}$.
\end{Theorem}

We sketch their proof of Theorem~\ref{ssthm} using the technology described here. In fact,
in \cite{HruSh334} they prove more -- If
$\{\Delta_m(\wbar_m):m\in\omega\}$ is any countable set of partial
types and there is an uncountable model $N$ of a countable,
superstable theory $T$ omitting each $\Delta_m$, then there is a
model $M$ of size continuum, again with an asymptotically similar
subset, that also omits each $\Delta_m$. As well, using the same
machinery they obtain the same conclusion for a countable stable
theory, at the cost of requiring the original model $N$ to have size
$\aleph_{\omega+1}$.

By employing the extensive calculus of non-forking, Shelah has gleaned many structural consequences
from his notion of a {\em stable system} of models.

\begin{Definition}  {\em Let $I$ be any non-empty index set.  A {\em stable system} of countable models of $T$ is a set
$\{M(s):s\in[I]^{<\omega}\}$ of countable models of $T$ satisfying:
\begin{itemize}
\item  If $s\subseteq t$, then $M(s)\preceq M(t)$;
\item  For all $s,t\in [I]^{<\omega}$, then $M(s)$ and $M(t)$ are independent (i.e., do not fork) over $M(s\cap t)$.
\end{itemize}
}
\end{Definition}

A primary tool for construction stable systems of models is {\em
domination}.  That is, given a pair of models $M\preceq M'$ and a
subset $B\subseteq M'$, we say {\em $B$ dominates $M'$ over $M$} if,
for any set $X$ (in some larger model), if $X$ is independent from
$B$ over $M$, then $X$ is independent from $M'$ over $M$. As we are
working over models in a stable theory, a sufficient condition for
domination is Lachlan's notion \cite{Lachlan2card} of locally atomic
models, $\ell$-atomicity:

\begin{Definition}  {\em Given a set $B$, a complete type $p\in S_n(B)$ is {\em locally ($\ell$-isolated)} if, for every partitioned formula $\phi(\xbar,\ybar)$, there is a formula $\psi(\xbar)\in p$ such that
$\psi(\xbar)\vdash \phi(\xbar,\bbar)$ for every $\phi(\xbar,\bbar)\in p$.  We call a model $M'$ {\em $\ell$-atomic over $B$} if, for every finite $\abar$ from $M'$, $\tp(\abar/B)$ is
$\ell$-isolated.
}
\end{Definition}

\noindent
A fundamental fact is that for stable theories,  if $M\subseteq B$ and if $M'$ is $\ell$-atomic over $B$, then $M'$ is dominated by $B$ over $M$.

Hrushovski and Shelah's proof of Theorem~\ref{ssthm} breaks into two pieces.
 The first part, which uses some highly technical stability-theoretic machinery
 (including the existence of definable groups in some instances) states that one can
 find a `very rich' stable system indexed by $I=\omega_1$ of elementary substructures of any uncountable model
 $N$ of a superstable theory $T$.


\begin{Theorem} \cite{HruSh334} \label{stablesystem} Let $N$ be an uncountable model of a countable, superstable theory $T$.
There is a stable system $\{M(s):s\in[\omega_1]^{<\omega}\}$ of countable, elementary substructures of $N$ and an independent subset $C=\{c_i:i\in\omega_1\}$ over $M(\emptyset)$
of $N$ that satisfy:
\begin{enumerate}
\item
For each $i\in\omega_1$, $c_i\in M(\{i\})$ and $M(\{i\})$ is $\ell$-atomic over $Mc_i$;
\item   For each $i\in\omega_1$ and
    $\theta(x,\bbar)\in\tp(c_i/M(\emptyset))$, there are infinitely many
    $j\in\omega_1$ such that $M(\{j\})\models\theta(c_j,\bbar)$;
 and
\item  For $|s|\ge 2$, $M(s)$ is $\ell$-atomic over
    $\bigcup\{M(t)\mcolon t\subsetneq s\}$.
\end{enumerate}
\end{Theorem}

%
%

As this theorem is rather technical, we only sketch the argument here and use some unexplained notation.

\medskip\noindent
{\bf Proof sketch.}  Without loss, we may assume $N$ has cardinality $\aleph_1$.
Fix an enumeration $\<a_i\mcolon i\in\omega_1\>$ of $N$.  For
each $i\in\omega_1$, let $A_i=\{a_j\mcolon j<i\}$ and let
$p_i=tp(a_i/A_i)$. As each $p_i$ is based on a finite set, for each
$i$ there is some $j<i$ such that $p_i$ is based on $A_j$.  By
Fodor's Lemma, there is some $j^*$ and a stationary subset
$S\subseteq\omega_1$ such that for each $i\in S$, $i>j^*$ and $p_i$ is
based on $A_{j^*}$.  Fix such a $j^*$ and put $B:=A_{j^*}$.
 So $B$ is
countable, and by reindexing $S$, we have an uncountable set
$C=\{c_i\mcolon i\in\omega_1\}$ that is independent over $B$.

Next, choose a countable $M\preceq N$ such that $B\subseteq M$ and $M$ is an $na$-substructure of $N$.  Using superstability, by removing at most countably many of the
$c_i$'s we obtain that the remaining, uncountably many elements are independent over $M$.

Now that we have chosen $M$ and $I$, it remains to construct our
stable system $\<M(s)\mcolon s\in[\omega_1]^{<\omega}\>$. But this
follows immediately by successive applications of the Corollary on page 302 of
\cite{HruSh334}.
\endproof

The second part of the proof of Theorem~\ref{ssthm} can be proved using the technology of this paper.
For this half, only stability is needed.
\begin{Theorem} \label{usestablesys} Suppose $T$ is a countable, stable theory and $\{M(s):s\in[\omega_1]^{<\omega}\}$ is a stable system of countable elementary submodels of an atomic model $N$
satisfying Clauses~(1)-(3) of Theorem~\ref{stablesystem}.  Then there is a Borel, atomic model $N_1$ of size continuum with an asymptotically similar subset
$\{a_\eta:\eta\in 2^\omega\}$.  More generally, if $N$ omits a countable set $\{\Delta_m:m\in\omega\}$ of types, then $N_1$ can be chosen to omit each $\Delta_m$.
\end{Theorem}

\bp For this
application, we take our set  $Z$  of variable symbols to be $X\cup Y$,
where $X=\{x_\eta:\eta\in 2^\omega\}$ and  $Y=\{y_{s,i}:s\in [2^\omega]^{<\omega}, i\in\omega\}$.

Choose any fmac $A\subseteq 2^{<\omega}$ with  an enumeration $\<a_j:j\in A\>$.  Suppose that $f:A\rightarrow\omega_1$ is any injective mapping.
Any such $f$ describes a finite tuple $\cbar_f:=\<c_{f(j)}:j\in A\>$ from the distinguished independent set $C=\{c_i:i\in\omega_1\}$.  Also, $f$ 
extends to
a map $f:\P(A)\rightarrow[\omega_1]^{<\aleph_0}$ by $f(t):=\{c_{f(j)}:j\in t\}$.

With this notation, define the set $\PP_A$ of $A$-commitments to be the set of instantiated $Z_A$-formulas $\phi(\xbar,\zbar)$, where $\zbar:=\<\ybar_s:s\subseteq A\>$
and each tuple $\ybar_s$ is from
$\{y_{s,i}:i\in\omega\}$, for which there is some injective $f:A\rightarrow\omega_1$ and tuples $\<\bbar_s:s\subseteq A\>$ from $M(f(s))$ so that
$N\models\phi(\cbar_f,\bbar_s:s\subseteq A)$.
As usual, let $(\PP,\le)$ be the partial order where $\PP=\bigcup\{\PP_A:A$ an fmac$\}$ and $\le$ is defined as in Section~\ref{fmac}.
As the given model $N$ and hence each of the submodels $M(s)$ omit each $\Delta_m$, the {\bf Omitting $\Delta_m$} conditions are easily verified.
As well, the verifications of the density conditions {\bf Completeness} and {\bf Henkin witnesses} are straightforward.  For both, fix an fmac $A$ and an $A$-commitment
$\phi(\xbar,\ybar)\in\PP_A$.  Choose an injective function $f:A\rightarrow \omega_1$ and tuples $\bbar_s$ from $M(f(s))$ such that $M(f(A))\models \phi(\cbar_f,\bbar_s:s\subseteq A)$.

\medskip
\noindent
{\bf Completeness:}  Choose any  instantiated $Z_A$-formula $\psi(\zbar)$ and partition its variables as $\psi(\xbar,\ybar_s:s\subseteq A)$.  By adding dummy variables to both $\phi$ and $\psi$, we may assume they have the same instantiated variables.  To decide how to extend $\phi$, we simply appeal to $M(f(A))$.  On one hand, if $M(f(A))\models\psi(\cbar_f,\bbar_s:s\subseteq A)$, then put $\phi^*:=\phi\wedge\psi$; put $\phi^*:=\phi\wedge\neg\psi$ otherwise.

\medskip
\noindent
{\bf Henkin witnesses:}  Choose any instantiated $Z_A$-formula $\theta(w,\zbar)$ with $w$ free.
In the notation of Definition~\ref{t(z)}, put $t:=t(\zbar)$.  Then the subsequence  $\dbar$ of $\<\cbar_f,\bbar_s:s\subseteq A\>$ corresponding to $\zbar$ is contained in $M(f(t))$.
 As above, there are two cases.  If $M(f(A))\models \neg \exists w \theta(w,\dbar)$,
then put $\phi^*:=\phi\wedge\neg\exists w \theta(w,\zbar)$.  Otherwise,
 append a new element $y_{t,j}$ to $\ybar_t$, forming $\ybar_t'$,
and put $\phi^*:=\phi\wedge\theta(y_{t,j},\zbar)$.  As $M(f(t))\preceq M(f(A))$, there is $b^*\in M(f(t))$ witnessing $\theta(w,\dbar)$.  This extra element witnesses that $\phi^*\in \PP_A$.

\medskip

By contrast, the verification of {\bf Splitting} is more involved, and requires new ideas.  As above, fix an enumerated fmac $A=\<a_i:i<n\>$ and an injective $f:A\rightarrow\omega_1$
that witnesses that $\phi(\xbar,\zbar)\in \PP_A$.  Choose an arbitrary $a\in A$, but to ease notation, suppose that $a=a_0$ and choose $\phi(\xbar,\zbar)\in\PP_A$.
As notation, let $A^-=A\setminus \{a_0\}$, let $A_0=A^-\cup\{a\conc 0\}$ and $A_1=A^-\cup\{a\conc 1\}$. Thus, $A^{*a}=A_0\cup A_1$
and the liftings $h_0,h_1:A\rightarrow A^{*a}$ map onto $A_0,A_1$, respectively. Fix an enumeration $\<s_i:i<2^n\>$ of $\P(A)$ that satisfies
(I) $i\le j$ whenever $s_i\subseteq s_j$ and (II)  the initial segment $\<s_i:i<2^{n-1}\>$ enumerates $\P(A^-)$.

Our first move is to `improve' our formula $\phi(\xbar,\zbar)\in\PP_A$.  As notation, for each $i<2^n$, let $\phi_i(\xbar,\ybar_j:j<i)$ be the restriction of $\phi$ to the smaller set of variables (we write $\ybar_j$ in place of the more cumbersome $\ybar_{s_j}$).  Call an $A$-commitment $\phi$ {\em self-sufficient} if, for every $0<i<2^n-1$,
$$\phi_i(\xbar,\ybar_j:j<i)\ \vdash\  \exists \ybar_i\, \phi_{i+1}(\xbar,\ybar_j:j\le i)$$
The notion of a self-sufficient commitment is a variant on what Hrushovski and Shelah call an `S-condition' in \cite{HruSh334}.  There, with Proposition~2.3(a) they prove:

\medskip

\noindent{\bf Claim:}  For any fmac $A$, every $\phi\in\PP_A$, has a self-sufficient $\phi^*\in\PP_A$ extending $\phi$.  Moreover, if $f:A\rightarrow\omega_1$ witnesses that
$\phi\in\PP_A$, then the same function $f$ witnesses that $\phi^*\in\PP_A$.

\medskip

Given the Claim, to verify {\bf Splitting} we may assume that $\phi$ itself is self-sufficient.  Choose an injective function $f:A\rightarrow\omega_1$ and tuples
$\bbar_i$ from $M(f(s_i)$ for each $i<2^n$ such that $N\models\phi(\cbar_f, \bbar_i:i<2^n)$, where $\bbar_i$ is short for $\bbar_{s_i}$.
 Getting half of the witnessing set is routine, and just amounts to adjusting the notation.  Let
$f_0:A_0\rightarrow \omega_1$ be defined as $f_0(a\conc 0)=f(a)$ and $f_0(a')=f(a')$ for all $a'\in A^-$.  In particular, $\cbar_{f_0}=\cbar_f$ so $f_0$ witnesses that
$h_0(\phi)$ is consistent.  Write $\cbar_f$ as $c_0\conc\cbar^*$.  The second half will require us to find an element $c'\in C\setminus \cbar_f$ so that $\tp(c'/M(\emptyset))$
is sufficiently close to $\tp(c_0/M(\emptyset))$ and then finding tuples $\<\bbar_i':2^{n-1}\le i< 2^n\>$ from the stable system.  First, note that $c_0$ is independent from $\cbar^*$ over $M(\emptyset)$.
Coupled with the fact that each $\bbar_i$ is dominated by $\{c_{f(a)}:a\in s_i\}$ over $M(\emptyset)$, there is a formula $\delta(x)\in\tp(c_0/M(\emptyset))$ so that if $c'$ is any realization of $\delta$ that is independent from $\cbar^*$ over $M(\emptyset)$, then
$$N\models \phi_i(c'\cbar^*,\bbar_j:j<i) \quad \hbox{for all $i<2^{n-1}$}$$
However, Clause~(1) of our hypotheses on our stable system imply that there is some $c_\beta\in C\setminus \cbar_f$ that satisfies these requirements.
Now, define $f_1:A_1\rightarrow \omega_1$ by $f_1(a\conc 1)=\beta$ and $f_1(a')=f(a')$ for all $a'\in A^-$.  Then, using the self-sufficiency of $\phi$, one recursively finds
tuples $\bbar_j'$ from $M(f_1(s_j))$  for each $2^{n-1}\le j < 2^n$ such that
$$N\models \phi_k(c_\beta\cbar^*,\<\bbar_i:i<2^{n-1}\>, \<\bbar'_j:2^{n-1}\le j < k\>) \ \hbox{for each $2^{n-1}\le k< 2^n$}$$
Combining these two halves yields that $f^*=f_0\cup f_1$ witnesses that $\phi':=h_0(\phi)\wedge h_1(\phi)\wedge x_{a\conc 0}\neq x_{a\conc 1}$ is in $\PP_{A^{*a}}$.

With the verification of {\bf Splitting} in hand, Theorem~\ref{usestablesys} and hence Theorem~\ref{ssthm} follow immediately by an application of Theorem~\ref{MAIN}.
\endproof

\begin{Remark}  {\em This result does not immediately translate to
the study of complete sentences of $L_{\omega_1,\omega}$. While
stability notions are defined in that context (\cite{Baldwincatmon}),
the superstability hypothesis on the ambient theory here is vastly
stronger than infinitary stability which concerns only the atomic models.
}

\end{Remark}

\section{Applications II -- Theories with Skolem functions}  \label{arbapp}

\numberwithin{Theorem}{subsection} \setcounter{Theorem}{0}


In this section we give applications of the Henkin method outlined in
the previous sections to construct customized models of size continuum of theories that have Skolem functions.  We first indicate how the existence of Skolem functions allows for a streamlining of our technique.  Recall that  if $T$ is a complete theory
 that has Skolem functions, then given any model $M$ of $T$, the Skolem hull of any subset $C\subseteq M$ will be an elementary substructure $N\preceq M$ in which each $b\in N$
is the interpretation of $\tau(c_1,\dots,c_k)$ for some $L$-term $\tau$ and some sequence $(c_1,\dots,c_k)$ of distinct elements of $C$.
In particular, having such tight control obviates the need for $Y$-variables!  More precisely, extra elements are needed to close $X$ to a model, but the existence of Skolem functions
makes their interpretations unique, and thus redundant.  Within this section, we will take $Z=X=\{x_\eta:\eta\in 2^\omega\}$ as our set of variable symbols
and we will construct a complete type $\Gamma(X)$ that is consistent with $T$.  As noted above, since $T$ admits
 Skolem functions,
 simply by taking the definable closure of any realization of $\Gamma$ inside any model, $\Gamma(X)$ uniquely determines a model of $T$.

 Thus, if $T$ has definable Skolem functions, then the {\bf Henkin witnesses} condition becomes vacuous.  As we are only concerned with $X$-variables,
the  {\bf Completeness} and {\bf Splitting} are easier to verify.
As usual, the {\bf Modeling $T$} clause is satisfied so long as every formula describing a commitment is satisfied in a model of $T$.
However, more care must be taken with {\bf Omitting $\Delta$}.  In particular,
our construction has to ensure that no $X$-instantiated  $L$-term $t(x_{\eta_1},\dots,x_{\eta_n})$ (or $m$-tuple of terms if $\Delta$ is $m$-ary) realizes $\Delta$.
In practice this will be easy to ensure, so long as the `witnessing models' each omit $\Delta$.

It might seem that  definable Skolem functions are in irreconcilable
conflict with the existence of large atomic models.  Indeed, if such a theory is countable, it cannot have an uncountable atomic model.
Despite that, we can use the technique here to construct atomic models of size continuum by expanding the language as follows.

\begin{Definition}  {\em  A {\em representation} of an
$L(\Phi)$-$L_{\omega_1,\omega}$-sentence  $\Phi$ is a triple $(L,T',\Delta(w))$ such that $L$ is a countable extension of $L(\Phi)$, $T$ is an $L$-theory,
and $\Delta(w)$ is a 1-type such that $Mod(\Phi)$ is equal to the class of $L(\Phi)$-reducts of models of $T'$ that omit $\Delta$.
Abusing notation somewhat, a {\em Skolemized representation} is a representation in which $T'$ admits definable Skolem functions,
admits elimination of quantifiers,
and has a pairing function.
}\end{Definition}

Applying Remark~\ref{easytrans}, it easy to find $(L,T',\Delta(w))$, a Skolemized representation, for an arbitrary complete $L_{\omega_1,\omega}$-sentence; choose a countable language
$L'\supseteq L$, a first-order $L'$-theory $T'$, and a partial  type $\Delta(w)$
such that the models of $\Phi$ are precisely the $L$-reducts of
models of $T'$ that omit $\Delta(w)$.  By expanding the language still
further (but maintaining countability)
 we may assume  $(L,T',\Delta(w))$ is a Skolemized representation.
 Then, if we construct a model $M'$ of $T'$  of size continuum
that omits each of the
partial types $\Delta_n$ given in the proof of Remark~\ref{smalltrans}, its reduct $M$ to $L$ is a large atomic model of $T$.

\subsection{Two-cardinal models}

In a pair of papers, \cite{Sh37,Sh49}, Shelah proves a celebrated
two-cardinal transfer theorem. For us, it is noteworthy as this is
apparently the first place where he uses the concept of asymptotic
similarity. In this situation we are able to simplify by assuming
Skolem functions as just discussed in the introduction.

Let $T$ be a theory in a countable language $L$ with a distinguished unary predicate $U$.
A model $M$ of $T$ is a {\em $(\kappa,\lambda)$-model} if $M$ has cardinality $\kappa$, but $|U(M)|=\lambda$.
We are interested in constructing a $(2^{\aleph_0},\aleph_0)$-model of $T$.  Clearly, we will not be able to succeed for an arbitrary theory $T$,
but we seek a sufficient condition on $T$  for a $(2^{\aleph_0},\aleph_0)$-model to exist.

Suppose that a countable theory $T$ has Skolem functions.  Thus, as suggested in the introduction to this section,
 take $X=\{x_\eta:\eta\in 2^\omega\}$ to be our distinguished set of variables, and let $\Gamma(X)$ be the
{\bf partial} type in these variables satisfying:

\begin{enumerate}
\item  $\neg U(x_\eta)$ and $x_\eta\neq x_{\eta'}$ for distinct $\eta,\eta'\in 2^\omega$;
\item  For each $k,\ell\in\omega$, for each $k$-ary $L$-term $\tau(w_1,\dots,w_k)$ and for all pairs of $\ell$-similar $k$-tuples $\etabar=(\eta_1,\dots,\eta_k)$ and
$\etabar'=(\eta_1',\dots,\eta_k')$ we have:
$$U(\tau(x_{\eta_1},\dots,x_{\eta_k}))\rightarrow \left[\tau(x_{\eta_1},\dots,x_{\eta_k})=\tau(x_{\eta_1'},\dots,x_{\eta_k'})\right]$$
\end{enumerate}
The following Lemma is immediate.

\begin{Lemma}[Shelah,\cite{Sh37}] \label{template}
Suppose that $T$ is a countable $L$-theory with Skolem functions.  If $T\cup\Gamma(X)$ is consistent,
then $T$ has a $(2^{\aleph_0},\aleph_0)$-model.
\end{Lemma}

\bp  Choose a model $M\models T$ with a subset $\{c_\eta:\eta\in 2^\omega\}$ satisfying $\Gamma(X)$.  It is easily
checked that  the Skolem hull of $\{c_\eta:\eta\in 2^\omega\}$ is a $(2^{\aleph_0},\aleph_0)$-model of $T$.
\endproof

But when is the type $\Gamma(X)$ consistent with $T$?    By compactness, it suffices to show that every finite subset of $\Gamma(X)$ is consistent with $T$.
That is, it suffices to show that every partial type $\Gamma_{\T}(X_F)$ is consistent with $T$, where {\em $\T$ is a finite set of $L$-terms }(of various arities), $F$ is a finite subset of $2^{\omega}$, and $\Gamma_{\T}(X_F)$ is the finite subset of $\Gamma(X)$ that mention only terms $\tau\in\T$ and variables $\{x_\eta:\eta\in F\}$.

For the remainder of this discussion,  fix a finite set $\T$ of $L$-terms.
Note that for any finite set $F\subseteq 2^\omega$, there is a $k<\omega$ such that $\{\eta|k:\eta\in F\}$ are distinct elements of $2^k$.
Choose any $m\ge k$, and consider the standard fmac $2^m\subseteq 2^{<\omega}$.
Let $\Gamma_{\T}(X_{2^m})$ be the set of $X_{2^m}$-instantiated formulas formed by replacing each variable symbol $x_\eta\in X_F$ by $x_{\eta|m}\in X_{(2^m)}$.
As any finite tuple $\cbar$ from any model $M\models T$ (indeed, any $L$-structure) realizes $\Gamma_{\T}(X_F)$ if and only if it realizes $\Gamma_{\T}(X_{2^m})$,
in order to show that $T\cup \Gamma_{\T}(X)$ is consistent, it suffices to prove that $T\cup \Gamma_{\T}(X_{2^m})$ is consistent for each of the standard fmacs $2^m$.

This overview of the proof was clear to Shelah at the time he wrote \cite{Sh37}, but it took him over a year to work out the combinatorics in \cite{Sh49} that led to the proof of
$(\aleph_\omega,\aleph_0)\rightarrow (2^{\aleph_0},\aleph_0)$.   We can now view his arguments as a slight variant on {\bf Splitting}.
Indeed, with our finite choice $\T$ of terms remaining fixed, choose any fmac $A$ (and an enumeration  $\<a_i:i<n\>$
thereof).  Suppose $M\models T$ and $\cbar=\<c_a:a\in A\>\in M^n$
is a tuple from $M$ realizing $\Gamma_{\T}(X_A)$.  Choose any $a\in A$ (say $a=a_j$).  We want to find an element $c^*\in M\setminus\{c_a:a\in A\}$ so that the
$(n+1)$-tuple $\cbar\conc c^*$ realizes $\Gamma_{\T}(X_{A^{*a}})$.  To obtain a sufficient condition for this, consider the equivalence relation $E_n$ on $(M)^n$,
the set of $n$-tuples of {\em distinct} elements from $M$ given by
$E_n(\cbar,\dbar)$ if and only if:
\begin{quotation}
\noindent For each $\tau(\wbar)\in\T$ and corresponding subsequences $\cbar',\dbar'$ with $\lg(\wbar)=\lg(\cbar')=\lg(\dbar')$,
 {\bf either} $M\models \neg U(\tau(\cbar'))\wedge \neg  U(\tau(\dbar'))$
 {\bf or}  $M\models \tau(\cbar')=\tau(\dbar')$.
\end{quotation}
It is easily verified that if $M$ is a $(\kappa,\lambda)$-model, then $E_n$ is an equivalence relation  on $(M)^n$ with at most $\lambda$ classes.
In terms of the discussion above, given $\cbar\in (M)^n$, we are seeking $c^*$ such that $E_n(\cbar,\cbar^*)$ holds, were $\cbar^*$ is formed by replacing $c_i$ by $c^*$ in
$\cbar$.

Finally, recall that every fmac $A$ can be constructed from $\{\<\>\}$ by a sequence of   $(|A|-1)$ splittings.  The following Proposition is merely a restatement of Theorem~5 of
\cite{Sh49}, noting that any equivalence relation $E_n$ on $(M)^n$ with at most $\lambda$ classes can be identified with a function $f\colon (M)^n\rightarrow\lambda$.

\begin{Proposition}[Shelah]  Fix any $m$, let $n=2^m-1$, and fix  a sequence $\<A_\ell:\ell\le n\>$ of fmacs  and a sequence $\<a_\ell:\ell<n\>$
such that
$A_0=\{\<\>\}$, $A_n=2^m$, and each $A_{\ell+1}=(A_\ell)^{*a_\ell}$.
 If $M$ is a $(\lambda^{+n},\lambda)$-model of $T$, then there is a tuple $\cbar=\<c_0,\cdots,c_n\>$ such that for every
$0<\ell\le n$,  $M\models E_\ell(\cbar\mr{\ell},\cbar\mr{\ell}^*)$, where $\cbar\mr{\ell}^*$ is obtained by substituting $c_n$ for the element of $\cbar_\ell$ coded by $a_\ell$.
In particular, for each $\ell\le n$, $\cbar_\ell$ realizes $\Gamma_{\T}(X_{A_\ell})$.
\end{Proposition}

Given this Proposition, the following Theorem of Shelah is immediate.

\begin{Theorem}[Shelah, \cite{Sh49}]   \label{2card2} $(\aleph_\omega,\aleph_0)\rightarrow (2^{\aleph_0},\aleph_0)$.  Indeed, if for every $n$, a theory $T$ admits a gap $n$-model,
i.e. a $((\lambda_n)^{+n},\lambda_n)$-model, then $T$ admits a $(2^{\aleph_0},\aleph_0)$-model.
\end{Theorem}

\bp  First, we may assume $T$ has Skolem functions.  Next, by Lemma~\ref{template} we need only show that $T\cup \Gamma(X)$ is consistent.
Fix any finite set $\T$ of $L$-terms.  By applying the Proposition for each $m$, we obtain the consistency of $T\cup\Gamma_{\T}(X_{2^m})$ for each of the standard
fmacs $2^m$, so we finish by compactness.
\endproof

%
%
%

The proof of Theorem~\ref{2card2} is an early exemplar of the `method
of identities' which has had many applications to prove two cardinal
theorems and compactness theorem in logics with generalized
quantifiers. See the account in \cite{ShVa}.

\subsection{What is the Hanf number for an atomic  model in the continuum?}
%
%
%
%
%

Classically, a `Hanf number' for a class of structures is the least cardinal $\lambda$ such that if the class of structures has one of size
$\lambda$, then it has arbitrarily large structures.  For example, Morley proved that if a sentence $\Phi$  of $L_{\omega_1,\omega}$ has a model of
size $\beth_{\omega_1}$, then $\Phi$ has arbitrarily large models.  Here, we vary the Hanf number question by asking for the smallest cardinal $\lambda$
for which the existence of a model of $\Phi$ of size $\lambda$ implies the existence of a model of size continuum.
Since every model of a complete $L_{\omega_1,\omega}$-sentence $\Phi$ is atomic (for a fixed expansion of the language of $\Phi$)
answering this question for a complete sentence  gives the Hanf number for atomic models in the continuum.

Clearly, the
value of $\lambda$ can vary, depending on the size of the continuum.
However, in \cite{Sh522}, Shelah defines (Definition~\ref{deflam}) a
cardinal $\lambda_{\omega_1}(\aleph_0)$ that is invariant under
c.c.c.\ forcings (hence by adding enough Cohen reals, we may assume
that $2^{\aleph_0}>\lambda_{\omega_1}(\aleph_0)$) and proves that if
a sentence $\Phi$ of $L_{\omega_1,\omega}$ has a model of size
$\lambda_{\omega_1}(\aleph_0)$, then it has a model of size
$2^{\aleph_0}$.

%


He defines what we  call (since it measures the ability to split in the sense here) a {\em splitting rank} for finite subsets of  $L$-structures $M$ in a countable language as follows:

%

\begin{Definition}  \label{sprk}  {\em For every non-empty, finite $B=\{b_0,\dots,b_{n-1}\}\subseteq M$, we define
the {\em splitting rank}, {\em sprk$(B,M)$}, by induction on $\alpha$ via the following clauses:

\begin{itemize}
\item  $\sprk(B,M)\ge 0$ if $B\cap\acl_M(\emptyset)=\emptyset$;
\item   For arbitrary $\alpha$, $sprk(B,M)\ge \alpha+1$ if and
    only if, for every $j<n$ and quantifier-free\footnote{The
    restriction to quantifier-free formulas is inessential in our
    applications here, but is stated in this manner to match the
    usage in \cite{Sh522}.}  $L$-formula
    $\phi(w_0,\dots,w_{n-1})$, there is $b^*_j\in (M\setminus B)$
    such that $$M\models \phi(b_0,\dots,
    b_j,\dots,b_{n-1})\leftrightarrow
    \phi(b_0,\dots,b^*_j,\dots,b_{n-1})$$ and
    $\sprk(Bb^*_j,M)\ge\alpha$; and
\item For $\alpha$ a non-zero limit, $\sprk(B,M)\ge\alpha$ if and only if $\sprk(B,M)\ge\beta$ for every $\beta<\alpha$.
\end{itemize}
}
\end{Definition}

Then define $\sprk(M)=\sup\{\sprk(B,M)+1:B\ \hbox{a finite subset of $M$}\}$ if the supremum exists,
or $\sprk(M)=\infty$ otherwise.

%

As extreme examples, suppose $B$ is a finite subset of $M$ satisfying $B\cap\acl(\emptyset)=\emptyset$, but some $b\in B$ is algebraic over $\bbar=B\setminus\{b\}$.
Then, if the formula $\phi(u,\bbar)$ witnesses the algebraicity, i.e., $M\models\phi(b,\bbar)\wedge\exists^{=k}u\phi(u,\bbar)$, then as successive splittings of this $B$ would require more and more distinct witnesses, we conclude that
$\sprk(B,M)<k$.  On the other extreme, an easy induction on $\alpha$ shows that $\sprk(B,M)\ge\alpha$ for any finite subset $B$ of any asymptotically similar subset
$\{a_\eta\colon \eta\in 2^\omega\}\subseteq M$ and any ordinal $\alpha$.   Thus, $\sprk(M)=\infty$ whenever $M$ contains an asymptotically similar subset.


%
%
%
%
%


\begin{Definition}\label{deflam} $\lambda_{\omega_1}(\aleph_0)$ is the least cardinal
$\lambda$ such that any structure $M$ of size
$\lambda$ for any countable language  necessarily has $\sprk(M)\ge\omega_1$.
\end{Definition}

 In \cite{Sh522},
Shelah proves that $\aleph_{\omega_1}\le
\lambda_{\omega_1}(\aleph_0)\le\beth_{\omega_1}$ and that this
cardinal is preserved under c.c.c.\ forcings. As the continuum can be
made arbitrarily large by adding enough Cohen reals (which is a
c.c.c.\ forcing) it is consistent that
$2^{\aleph_0}>\lambda_{\omega_1}(\aleph_0)$. Despite considerable
work on the problem, the question
\begin{quotation}
`Does $ZFC$ prove that $\lambda_{\omega_1}(\aleph_0)=\aleph_{\omega_1}$?'
\end{quotation}
remains open.
He also gives examples of sentences $\Phi_\alpha$ of $L_{\omega_1,\omega}$ for each $\alpha<\omega_1$ such that each $\Phi_\alpha$ has a model $M$ with $\sprk(M)=\alpha$
and no models of larger splitting rank; thus, in general, $\lambda_{\omega_1}(\aleph_0)\ge\aleph_{\omega_1}$.
The main theorem of \cite{Sh522} is a pleasant application of the methods developed
 in the previous sections:



\begin{Theorem}[Shelah,\cite{Sh522}] \label{Hanf}
has a model $M$ of size at least  $\lambda_{\omega_1}(\aleph_0)$, then
  $\Phi$ has a Borel model of size continuum
that contains an asymptotically similar subset $\{c_\eta:\eta\in 2^\omega\}$.
\end{Theorem}

\bp
Let  $(L,T',\Delta(w))$ be a Skolemized representation of $\Phi$.
%
As $T'$ has Skolem functions, take
 $Z=X=\{x_\eta:\eta\in 2^\omega\}$.  We will construct a complete type $\Gamma(Z)$ that is consistent with $T'$ and such that,
if $N$ is any model of $T'$ and $\{c_\eta:\eta\in 2^\omega\}$ realizes $\Gamma(Z)$ in $N$, then the Skolem hull of $\{c_\eta:\eta\in 2^\omega\}$ will omit $\Delta(w)$.


To accomplish this, for each fmac $A$ of $2^{<\omega}$, let $\PP_A$ denote all instantiated formulas $\phi(\xbar)=\phi(x_a:a\in A)$ that satisfy:
\begin{quotation}
\noindent For every $\alpha<\omega_1$ there is some  $\bbar_\alpha$ from $M'$ that realizes $\phi$ and
such that \hbox{$\sprk(M',\bbar_\alpha)\ge\alpha$.}
\end{quotation}
Take $\PP=\bigcup\{\PP_A:A\ \hbox{an fmac}\}$ and define  $\le$ to be the usual extension relation on commitments given in Section~4.

  As $M'$ is a model of $T'$, the structures we build will be models of $T'$.
Also, as   $T'$ has Skolem functions, the {\bf Henkin witnesses} conditions are trivial.  More interesting verifications are:

\medskip\noindent {\bf Completeness:}  Fix an fmac $A$ and an $A$-commitment $\phi(x_a:a\in A)\in\PP_A$, and choose any instantiated $X_A$-formula $\psi(x_a:a\in A)$.
As $\phi\in \PP_A$, for each $\alpha<\omega_1$, choose  $\bbar_\alpha$ from $M'$ realizing $\phi(\xbar)$
with $\sprk(M',\bbar_\alpha)\ge\alpha$.  There are now two cases:
First, if $Y=\{\alpha<\omega_1:M'\models\psi(\bbar_\alpha)\}$ is uncountable, then put $\phi^*(\xbar):=\phi\wedge\psi$.  By passing to this uncountable collection, it is evident that $\phi^*\in\PP_A$.  On the other hand, if $Y$ is countable, then as its complement is uncountable, put $\phi^*(\xbar):=\phi\wedge\neg\psi$ and again, $\phi^*\in\PP_A$
and extends $\phi$.

\medskip
The verification of {\bf Omitting $\Delta$} is similar.

\medskip\noindent{\bf Omitting $\Delta$:}  Given an fmac $A$ and $\phi\in\PP_A$, choose any $X_A$-instantiated $L$-term
$t(x_a:a\in A)$.  As above, for each $\alpha<\omega_1$ choose  a realization $\bbar_\alpha$
of $\phi$ in $M'$ with $\sprk(M,\bbar_\alpha)\ge\alpha$.  As $M'$ omits $\Delta(w)$, for every $\alpha$
there is $\delta_\alpha(w)\in\Delta$ such that $M'\models\neg\delta_\alpha(t(\bbar_\alpha))$.  As $\Delta$ is countable, choose a single $\delta^*\in\Delta$
such that $\{\alpha<\omega_1:M'\models\neg\delta^*(t(\bbar_\alpha))\}$ is uncountable.  Put $\phi^*(\xbar):=\phi\wedge\neg\delta^*(t(\xbar))$, which clearly extends
$\phi(\xbar)$.  By reindexing, it is evident that $\phi^*\in\PP_A$.

\medskip

The `shift' that occurs in the verification of {\bf Splitting} is reminiscent of the proof of Morley's Omitting Types theorem.

\medskip\noindent{\bf Splitting:}  Fix any fmac $A$, any $A$-commitment $\phi(\xbar)$, and choose any $a\in A$.
As in Definition~\ref{splittingdef}, let $A^{*a}=A\setminus\{a\}\cup\{a\conc 0,a\conc 1\}$, and put
$\phi^*:=\phi(h_0(\xbar))\wedge\phi(h_1(\xbar))\wedge x_{\delta\conc 0}\neq x_{\delta\conc 1}$.  It suffices to show that $\phi^*\in\PP_{A^{*a}}$.
To see this, for each $\alpha<\omega_1$,
choose  $\bbar_\alpha$ such that $M'\models\phi(\bbar_\alpha)$ and $\sprk(M',\bbar_\alpha)\ge\alpha+1$.
As $T'$ admits elimination of quantifiers, it follows from the definition of $\sprk$ that there is a 1-point extension $\bbar'_\alpha$ from $M'$ extending
$\bbar_\alpha$ that realizes $\phi^*$ with $\sprk(M',\bbar'_\alpha)\ge\alpha$.  Thus, $\phi^*\in\PP_{A^{*a}}$.

\bigskip
Once all of these conditions are satisfied, it follows from Theorem~\ref{MAIN} that there is a  Borel model $N^*$ of size continuum that models $T'$ and omits $\Delta(w)$
with an asymptotically similar subset
$\{c_\eta:\eta\in 2^\omega\}$.  As $T'$ has 
Skolem functions, the substructure $N'\preceq N^*$ generated by $\{c_\eta:\eta\in 2^\omega\}$ also models $T'$ and omits
$\Delta(w)$.  Thus, as explained in the introduction to Section~\ref{arbapp}
the reduct $N$  of $N'$ to the original language $L$ is a  Borel model of $\Phi$ that has both size continuum and an asymptotically similar subset.
\endproof

In \cite{Sh522}, Shelah draws an immediate Corollary from Theorem~\ref{Hanf}.  Given what we have proved above, all that is required is to code the hypotheses into a suitable
structure of cardinality $\lambda_{\omega_1}(\aleph_0)$.

\begin{Corollary}[Shelah]  Let $B\subseteq 2^\omega\times 2^\omega$ be a Borel subset of the product.  If $B$ contains a $\lambda_{\omega_1}(\aleph_0)$-square
(i.e., a subset $E\subseteq 2^\omega$ of size $\lambda_{\omega_1}(\aleph_0)$ such that
 $E\times E\subseteq B$) then there is a perfect subset $E^*$ of the continuum with $E^*\times E^*\subseteq
 B$.
\end{Corollary}

Recall that classically, Morley's Omitting Types theorem states that if there is a model of power $\beth_{\omega_1}$ omitting a type, then there are Ehrenfeucht-Mostowski models
that also omit the type.  However, by looking more closely at the proof, the hypotheses can be weakened to:  `For every $\alpha<\omega_1$, there is a model $M_\alpha$
of power at least $\beth_\alpha$ that omits the type.'
We note a similar analogy gives the following strengthening of
 Theorem~\ref{Hanf}.  Specifically, to prove Theorem~\ref{better2}, take, for each fmac $A$,  $\PP_A$ to be the set of all formulas $\phi(x_a:a\in A)$ such that for each $\alpha<\omega_1$, there is $\beta(\alpha)\ge\alpha$
and  $\bbar_\alpha$ from $M_{\beta(\alpha)}$ realizing $\phi$.

\begin{Theorem}  \label{better2} Suppose a sentence $\Phi$ of $L_{\omega_1,\omega}$ has a Skolemized representation $(L,T',\Delta)$. 
If, for every $\alpha<\omega_1$ there is a model $M_\alpha$ of $T'$ that omits $\Delta(w)$ and has $\sprk(M_\alpha)\ge\alpha$,
then there is a model $N$ of $T$ of size continuum that omits $\Delta(w)$ and has an asymptotically similar subset $\{c_\eta:\eta\in 2^\omega\}$.
\end{Theorem}

Theorem~\ref{better2} entails the following amusing Corollary.

\begin{Corollary} \label{amusing} Let $\Phi$ be any sentence of $L_{\omega_1,\omega}$ with  a Skolemized representation $(L,T',\Delta)$.
  If there is a $ZFC$-proof of the existence of a model of $\Phi$ of size continuum, then
there is a  Borel model of $\Phi$ with an asymptotically similar subset $\{c_\eta:\eta\in 2^\omega\}$.
\end{Corollary}

\bp  
It is easily seen by induction on $\alpha$ that for every $\alpha<\omega_1$ there is an $L_{\omega_1,\omega}$-sentence $\Psi_\alpha$ in the language $L'$
such that
an $L'$-structure $N'\models\Psi_\alpha$ if and only if $\sprk(N')\ge\alpha$.

To begin the proof of the Corollary, by forcing enough Cohen reals, work in a model $\V[G]$ of $ZFC$ in which
$2^{\aleph_0}>\lambda_{\omega_1}(\aleph_0)$.  As our forcing has the c.c.c., $\omega_1^{\V[G]}=\omega_1$.
Working in $\V[G]$, choose a model $N\models\Phi$ of size continuum.  Let $N'\models T'$ be an expansion of $N$ to $L'$ that omits $\Delta$.

As $|N'|\ge\lambda_{\omega_1}(\aleph_0)$,
$N'\models\Psi_\alpha$ for every $\alpha<\omega_1$.  Thus, for each $\alpha$, $\Phi\wedge\Psi_\alpha$ is formally consistent.  So, working in $\V$, for each $\alpha<\omega_1$
an application of  Karp's Completeness Theorem yields a (countable) model $M'_\alpha\models\Phi\wedge\Psi_\alpha$.
Collectively,  expansions of the models $\{M'_\alpha:\alpha<\omega_1\}$ satisfy the hypotheses of Theorem~\ref{better2}, so we finish.
\endproof


\begin{Remark}  {\em Both Theorem~\ref{Hanf} and Corollary~\ref{amusing} have analogues for atomic models.  Indeed, given a countable, complete theory $T$,
let $T'$ be a Skolemization of $T$ and let $\{\Delta_n\}$ be the partial types given at the end of  Subsection~\ref{atomic}.  Let $\Phi$ be the sentence of $L_{\omega_1,\omega}$
given in Remark~\ref{smalltrans} (with respect to $T'$).  Then, the $L$-reduct of any model $M'$ of $\Phi$ will be an atomic model of $T$; and conversely, every atomic model $M$
of $T$ has an expansion to a model $M'$ of $\Phi$.  Thus, it follows from Theorem~\ref{Hanf} that
 if a countable, complete, first order theory $T$ has an atomic model of size
$\lambda_{\omega_1}(\aleph_0)$, then $T$ has a Borel atomic model of size continuum.  Similarly, the analogue of Corollary~\ref{amusing} is that if there is a $ZFC$ proof of the existence of an atomic model  of size continuum for a  countable, complete, first order $T$, then there is a Borel, atomic model of $T$ of size continuum with an asymptotically similar subset.
}
\end{Remark}


\begin{Remark}  {\em
A glance at the definitions  shows that having
definable Skolem functions  is the antithesis of $\dcl$-triviality  (see Section~\ref{tttheory}).    In fact, the lack of non-trivial algebraic formulas
directly implies that every finite subset of $M$ has unbounded
splitting rank, i.e., $\sprk(A,M)=\infty$
for every finite subset $A$ of $M$.  In fact, this `arbitrary
splitting; condition characterizes  trivial-$\dcl$.  In fact, we have two proofs that theories with trivial $\dcl$ have atomic models in the continuum.  The first (Subsection~\ref{tttheory}) took place in a extension of the given vocabulary by predicates {\em definable in $L_{\omega_1,\omega}$}.  But the result also follows from the methods of this section using the next easy Fact and the fact that uncountable splitting rank gives a model in the continuum.


\begin{Fact}  \label{manydcl}
The following are equivalent for an $L$-structure $M$:
\begin{enumerate}
\item  $M$ has trivial $\dcl$;
\item  $\acl(A)=A$ for all subsets $A\subseteq M$;
\item For every finite subset $A\subseteq M$, $\sprk(A,M)\ge 1$;
\item  For every finite subset $A\subseteq M$, $\sprk(A,M)=\infty$;
\item  $(M,=)$ is a sufficient pregeometry.
\end{enumerate}
\end{Fact}
}
\end{Remark}

\end{document}